\documentclass{birkmult}
\usepackage{mathrsfs,amssymb}
 \newtheorem{thm}{Theorem}[section]
 
 \newtheorem{lem}[thm]{Lemma}
 \newtheorem{prop}[thm]{Proposition}
 \newtheorem{pbm}{Open Problem}
 \newtheorem*{taylor}{Taylor's Theorem}
 \newtheorem*{frostman}{Frostman's Lemma}
 
 \theoremstyle{definition}{
 \newtheorem{defn}[thm]{Definition}}
 
 \theoremstyle{remark}{
	 \newtheorem{rem}[thm]{Remark}
	 \newtheorem{rems}[thm]{Remarks}
	 }

 \numberwithin{equation}{section}
 
 \newcommand{\e}{\epsilon}
 \newcommand{\s}{\sigma}
 \newcommand{\g}{\mathbf{g}}
 \newcommand{\R}{\mathbf{R}}
 \newcommand{\E}{\mathrm{E}}
 
 \newcommand{\F}{\mathcal{F}}
 
 \newcommand{\M}{\mathrm{M}}
 \newcommand{\K}{\mathrm{K}}
 \newcommand{\dimh}{\dim_{_{\rm H}}}
 
 \renewcommand{\P}{\mathrm{P}}

\begin{document}

\title[Slices of Brownian Sheet]{Slices of Brownian Sheet:\\
	New Results, and Open Problems}
\author[D. Khoshnevisan]{Davar Khoshnevisan}
\address{%
	Department of Mathematics\\
	The University of Utah\\
	155 S 1400 E Salt Lake City, UT 84112--0090, 
	USA}
\email{davar@math.utah.edu}
\urladdr{http://www.math.utah.edu/\~{}davar}

\thanks{Research supported in part by a grant from the
	United States National Science Foundation}
\date{October 24, 2005}
%\dedicatory{This paper is dedicated to a memory of a memory}

\begin{abstract}
	We can view Brownian sheet as a sequence of interacting
	Brownian motions or \emph{slices}. Here we present
	a number of results about the slices of the sheet.
	A common feature of our results is that they exhibit
	phase transition. In addition,
	a number of open problems are presented.
\end{abstract}
\maketitle
\tableofcontents

\section{Introduction}

Let $B:=\{B(s\,,t)\}_{s,t\ge 0}$ denote a two-parameter Brownian sheet
in $\R^d$. That is, $B$ is a centered Gaussian process with
covariance matrix,
\begin{equation}
	\mathrm{Cov}\left( B_i(s\,,t) ~,~ B_j(u\,,v) \right)
	= \min(s\,,u)\cdot\min(t\,,v)\cdot \delta_{i,j}.
\end{equation}
We can assume without loss of generality that $B$ is  continuous.
Moreover, it is convenient to think of $B$ as the
distribution function of a $d$-dimensional white noise $\hat{B}$
on $\R^2_+$; i.e., we may think of $B(s\,,t)$ as
\begin{equation}
	B(s\,,t) = \hat{B} \left( [0\,,s]\times[0\,,t] \right).
\end{equation}
These properties were discovered 
first in \v{C}entsov \cite{Centsov}.

Choose and fix some number $s>0$.
The \emph{slice of $B$ along $s$} is
the stochastic process $\{ B(s\,,t)\}_{t\ge 0}$.
It is easy to see that if $s$ is non-random
then the slice of $B$ along $s$ is a scaled
Brownian motion. More precisely,
$t\mapsto s^{-1/2} B(s\,,t)$ is standard $d$-dimensional
Brownian motion. It is not  too difficult to
see that  if $s$ is random,
then the slice along $s$ need not be a Brownian 
motion. For instance,
the slice along a non-random $s$ hits points if and only if 
$d=1$. But there are random values of $s$ such that the
slice along $s$ hits zero up to dimension $d=3$; see
\eqref{eq:Orey--Pruitt} below.
Nonetheless, one may expect
the slice along $s$ to look like Brownian motion
in some sense, even for some random values of $s$.
[For example, all slices share the Brownian property that
they are continuous paths.]

A common question in infinite-dimensional
stochastic analysis is to ask if there are slices that behave
differently from $d$-dimensional Brownian motion
in a predescribed manner. There is a large literature
on this subject; see the survey
paper \cite{Kh04}. In this paper we present some
new examples where there is, generally, a ``cut-off
phenomenon'' or ``phase transition.''

Our first example is related to the zero-set
of the Brownian sheet. Orey and Pruitt \cite{OreyPruitt}
have proven that $B^{-1}\{0\}$ is non-trivial if and only if
the spatial dimension $d$ is three or less.
That is,
\begin{equation}\label{eq:Orey--Pruitt}
	\P\left\{ B(s\,,t)=0 \text{ \rm for some }
	s,t>0 \right\}>0\ \Longleftrightarrow\ 
	d\le 3.
\end{equation}
See also
Fukushima \cite{Fukushima} and Penrose \cite{Penrose}.
Khoshnevisan \cite{Kh99} has derived the following
refinement: For all non-random,
compact sets $E,F\subset(0\,,\infty)$,
\begin{equation}\label{eq:Kh:Cap-d/2}
	\P\left\{ B^{-1}\{0\} \cap
	(E\times F) \neq\varnothing
	\right\}>0\ \Longleftrightarrow\
	\text{Cap}_{d/2}(E\times F)>0,
\end{equation}
where $\text{Cap}_\beta$ denotes
``$\beta$-dimensional Riesz capacity.''
[These capacities are recalled in 
the appendix.]
The Orey--Pruitt theorem \eqref{eq:Orey--Pruitt}
follows immediately from \eqref{eq:Kh:Cap-d/2}
and Taylor's theorem [Appendix \ref{app:cap}].

Now consider the projection $\mathcal{Z}_d$ of
$B^{-1}\{0\}$ onto the $x$-axis. That is,
\begin{equation}
	\mathcal{Z}_d := \left\{
	s \ge 0 :\,
	B(s\,,t)=0\text{ \rm for some }t > 0 \right\}.
\end{equation}
Thus, $s\in \mathcal{Z}_d$ if and only if
the slice of $B$ along $s$ hits zero. Of course,
zero is always in $\mathcal{Z}_d$, and the latter is
a.s.\ closed. Our first result
characterizes the polar sets of $\mathcal{Z}_d$.

\begin{thm}\label{thm:BSMain-Intro}
	For all non-random, compact sets $F\subset (0\,,\infty)$,
	\begin{equation}
	\P\left\{  \mathcal{Z}_d  \cap F  \neq\varnothing \right\}>0
	\ \Longleftrightarrow\
	\text{\rm Cap}_{(d-2)/2} (F)>0.
	\end{equation}
\end{thm}

Theorem \ref{thm:BSMain-Intro} and Taylor's theorem
[Appendix \ref{app:cap}] together provide us with a new
proof of the Orey--Pruitt theorem \eqref{eq:Orey--Pruitt}.
Furthermore, we can apply a
codimension argument \cite[Theorem 4.7.1, p.\ 436]{Kh}
to find that
\begin{equation}\label{eq:dimh-Intro}
	\dimh \mathcal{Z}_d = 1 \wedge
	\left(2-\frac d2\right)^+\qquad\text{a.s.},
\end{equation}
where $\dimh$ denotes Hausdorff dimension
[Appendix \ref{app:dimh}].
Consequently, when $d\in\{2\,,3\}$,
the [Hausdorff] dimension of $\mathcal{Z}_d$
is equal to $2-(d/2)$. Oddly enough, this
is precisely the dimension of $B^{-1}\{0\}$ as well;
see Rosen \cite{Rosen,Rosen83}. But $\mathcal{Z}_d$ is
the projection of $B^{-1} \{0\}$ onto the $x$-axis.
Therefore, one might guess that $B^{-1}\{0\}$ and
$\mathcal{Z}_d$ have the
same dimension because all slices of $B$ have
the property that their zero-sets have zero dimension.
If $B$ were a generic function of two variables, then
such a result would be false, as there are simple
counter-examples. Nevertheless, the ``homoegenity''
of the slices of $B$ guarantees that our
intuition is correct in this case.

\begin{thm}\label{thm:trace}
	If $d\in\{2\,,3\}$, then the following holds
	outside a single $\P$-null set:
	\begin{equation}\label{eq:trace}
		\dimh \left( B^{-1}\{0\}
		\cap \left( \{s\}\times(0\,,\infty)\right) \right) =0
		\qquad\text{for all }s>0.
	\end{equation}
\end{thm}

\begin{rems}\begin{enumerate}
	\item Equation \eqref{eq:trace} is not valid when $d=1$.
		In that case, Penrose \cite{Penrose} proved that
		$\dimh ( B^{-1}\{0\}\cap(\{s\}\times(0\,,\infty)))
		=1/2$ for all $s>0$. In particular, Penrose's
		theorem implies that 
		$\mathcal{Z}_1=\R_+$ a.s.; the latter follows
		also from an earlier theorem of
		Shigekawa \cite{Shigekawa}.
	\item Almost surely, $\mathcal{Z}_d=\{0\}$
		when $d\ge 4$; see \eqref{eq:Orey--Pruitt}.
		This and the previous remark together show that
		``$d\in\{2\,,3\}$'' covers the only interesting
		dimensions.
	\item The fact that
		Brownian motion misses singletons in $\R^d$ for
		$d\ge 2$ implies that the Lebesgue measure
		of $\mathcal{Z}_d$ is a.s.\ zero
		when $d\in\{2\,,3\}$.
	\item It is not hard to see that the probability in Theorem
		\ref{thm:BSMain-Intro} is zero or one. 
		Used in conjunction
		with Theorem \ref{thm:BSMain-Intro}, this observation
		demonstrates that $\mathcal{Z}_d$ is a.s.\
		everywhere-dense when $d\le 3$.
\end{enumerate}\end{rems}

Next, we consider the random set,
\begin{equation}
	\mathcal{D}_d := \left\{ 
	s\ge 0 :\,  B(s\,,t_1)=B(s\,,t_2) \text{ for some }
	t_2>t_1>0\right\}.
\end{equation}
We can note that $s\in\mathcal{D}_d$
if and only if the slice of $B$ along $s$
has a double point.
 
Lyons \cite{Lyons} has proven that 
$\mathcal{D}_d$ is non-trivial if and only if
$d\le 5$. That is,
\begin{equation}\label{eq:lyons}\begin{split}
	\P\left\{ \mathcal{D}_d\neq\{0\}
	\right\}>0\ \Longleftrightarrow\
	d\le 5.
\end{split}\end{equation}
See also Mountford \cite{Mountford}.
Lyons's theorem \eqref{eq:lyons} is an improvement
to an earlier theorem
of Fukushima \cite{Fukushima}
which asserts the necessity of the condition ``$d\le 6$.''
Our next result characterizes the polar sets of
$\mathcal{D}_d$.

\begin{thm}\label{thm:Double-Intro}
	For all non-random, compact sets $F\subset (0\,,\infty)$,
	\begin{equation}
		\P\left\{ \mathcal{D}_d \cap F
		\neq\varnothing  \right\}>0\ \Longleftrightarrow\
		\text{\rm Cap}_{(d-4)/2} (F)>0.
	\end{equation}
\end{thm}
Lyons's theorem \eqref{eq:lyons} follows at once
from this and Taylor's theorem.  In addition, a codimension
argument reveals that  almost surely,
\begin{equation}\label{eq:dimD}
	\dimh \mathcal{D}_d = 1
	\wedge \left( 3 -\frac d2\right)^+.
\end{equation}
This was derived earlier by Mountford \cite{Mountford}
who used different methods.

\begin{rem}
	Penrose \cite{Penrose,Penrose89} has
	shown that $\mathcal{D}_d=\R^d_+$ a.s.\ when $d\le 3$.
	Also recall Lyons' theorem \eqref{eq:lyons}. Thus,
	Theorem \ref{thm:Double-Intro} has content only 
	when $d\in\{4\,,5\}$.
\end{rem}

In summary, our
Theorems \ref{thm:BSMain-Intro} and 
\ref{thm:Double-Intro} state that certain
unusual slices of the sheet can be found in 
the ``target set'' $F$ if and only if
$F$ is sufficiently large in the sense of capacity.
Next we introduce a property which is related to
more delicate features of the set $F$. Before 
doing so, let us set $d\ge 3$
and define
\begin{equation}\label{eq:rate_c}
	\mathcal{R}(s) := \inf\left\{
	\alpha>0:\
	\liminf_{t\to \infty} \frac{(\log t)^{1/\alpha}}{t^{1/2}}
	|B(s\,,t)| <\infty\right\}
	\qquad\text{for all }s>0.
\end{equation}
Thus, $\mathcal{R}(s)$ is the critical
escape-rate---at the logarithmic level---for the
slice of $B$ along $s$. Because
$t\mapsto s^{-1/2} B(s\,,t)$
is  standard Brownian motion for all fixed $s>0$,
the integral test 
of Dvoretzky and Erd\H{o}s \cite{DvoretzkyErdos} 
implies that
\begin{equation}\label{eq:DE}
	\P\left\{ \mathcal{R}(s) = d-2\right\}=1
	\quad\text{for all }s>0.
\end{equation}
That is, the typical slice of $B$ escapes at log-rate
$(d-2)$. This leads to the question,
``When are all slices of $B$ transient''?
Stated succinctly, the answer is: ``If and only if
$d\ge 5$.'' See Fukushima \cite{Fukushima} for
the sufficiency of the condition ``$d\ge 5$,''
and K\^ono \cite{Kono} for the necessity.
Further information can be found in
Dalang and Khoshnevisan \cite{DalangKh}.
Next we try to shed further light on
the rate of convergence of the transient slices of $B$.
Our characterization is in terms
of packing dimension $\dim_{_{\rm P}}$, which is recalled
in Appendix \ref{subsec:dimp}.

\begin{thm}\label{thm:DE}
	Choose and fix $d\ge 3$,
	and a non-random compact set $F\subset(0\,,\infty)$.
	Then with probability one:
	\begin{enumerate}
		\item $\mathcal{R}(s) \ge d-2-2\dim_{_{\rm P}}F$
			for all $s\in F$.
		\item If $\dim_{_{\rm P}} F <(d-2)/2$, 
			then $\mathcal{R}(s) = d-2-2\dim_{_{\rm P}}F$
			for some $s\in F$.
	\end{enumerate}
\end{thm}

\begin{rem}
	The condition that $\dim_{_{\rm P}} F<(d-2)/2$
	is always met when $d\ge 5$.
\end{rem}

The organization of this paper is as follows:
After introducing some basic real-variable
computations in Section \ref{sec:RV},
we prove Theorem \ref{thm:BSMain-Intro}
in Section \ref{sec:BSMain}. Our argument
is entirely harmonic-analytic and does not
require any probability theory; this proof rests
on a projection theorem for capacities
which may be of independent interest.
Theorems \ref{thm:Double-Intro} and \ref{thm:trace} are
respectively proved in  Sections \ref{sec:Double} 
and \ref{sec:trace}. Section \ref{sec:moreDouble}
contains a variant of Theorem \ref{thm:Double-Intro},
and Section \ref{sec:rates}
contains the proof of Theorem \ref{thm:DE}
and much more.
There is also a final Section \ref{sec:open}
wherein we  record some open problems.

Throughout, any $n$-vector $x$ is written, coordinatewise,
as $x= (x_1,\ldots,x_n)$. Moreover, $|x|$ will always
denote the $\ell^1$-norm of $x\in\R^n$; i.e.,
\begin{equation}
	|x|:=|x_1|+\cdots+|x_n|.
\end{equation}

Generic constants that do not depend on anything
interesting are denoted by $c,c_1,c_2,\ldots$; they are
\emph{always} assumed to be positive and finite, and their
values may change between, as well as within, lines.

Let $A$ denote a Borel set in
$\R^n$. The collection of all Borel probability measures on 
$A$ is always denoted by $\mathcal{P}(A)$.\\

\noindent\textbf{Acknowledgement.} A large portion of
this work was
motivated by enlightening discussions
with Robert Dalang over a period of several years.
A great many thanks are due to him.

\section{Preliminary Real-Variable Estimates}\label{sec:RV}

Our analysis depends on the properties
of three classes of functions. We develop the
requisite estimates here in this section.
Aspects of these lemmas overlap with Lemmas
1.2 and 2.5
of Dalang and Khoshnevisan \cite{DalangKh}.

Here and throughout, we define for all
$\e>0$ and $x\in\R$,
\begin{equation}\label{eq:f}\begin{split}
	f_\e(x) &:=  \left( \frac{\e}{
		|x|^{1/2}} \wedge 1\right)^d,\\
	F_\e(x) &:= \int_0^1 f_\e(y+|x|)\, dy,\\
	G_\e(x) &:= \int_0^1 F_\e(y+|x|)\, dy.
\end{split}\end{equation}

Our first technical lemma attaches a ``meaning''
to $f_\e$.

\begin{lem}\label{lem:g-f}
	Let $\g$ denote a $d$-vector of i.i.d.\
	standard-normal variables. Then there exist a constant $c$
	such that for all $\s,\e>0$,
	\begin{equation}
		c f_\e(\s^2)\le 
		\P\left\{ \s|\g|\le \e \right\}
		\le f_\e(\s^2).
	\end{equation}
\end{lem}

\begin{proof}
	This is truly an elementary result. However,
	we include a proof to acquaint the reader
	with some of the methods that we use later on.
	
	Let $M:= \max_{1\le i\le d}|\g_i|$,
	and note that  $|\g|\ge M$.
	Therefore,
	\begin{equation}
		\P\left\{ \s|\g|\le \e\right\} \le
		\left(\int_{-\e/\s}^{\e/ \s } 
		\frac{e^{-u^2/2}}{(2\pi)^{1/2}} \, du
		\right)^d\le \left(\frac{\e}{\s}\right)^d,
	\end{equation}
	because $(2/\pi)^{1/2}\exp(-u^2/2)\le 1$.
	The upper bound of the lemma
	follows because $\P\{\s|\g|\le\e\}$ is 
	also at most one. 
	To derive the lower bound we use the
	inequality $|\g|\le Md$ to find that
	when $\e\le \s$,
	\begin{equation}\begin{split}
		\P\left\{ \s|\g|\le \e\right\} &\ge
			\left(\int_{-\e/(\s d)}^{\e/(\s d)} 
			\frac{e^{-u^2/2}}{(2\pi)^{1/2}} \, du
			\right)^d
			\ge \left(\frac{2}{\pi d^2}\right)^{d/2} e^{-1/(2d^2)}
			\left( \frac{\e}{\s} \right)^d\\
		&= \left(\frac{2}{\pi d^2}\right)^{d/2} e^{-1/(2d^2)}
			f_\e(\s^2) := c_1 f_\e(\s^2).
	\end{split}\end{equation}
	The same reasoning shows that when $\e>\s$,
	\begin{equation}\begin{split}
		\P\left\{ \s|\g|\le \e\right\} &\ge
			\left(\int_{-1}^1 \frac{e^{-u^2/2}}{(2\pi)^{1/2}}
			\, du\right)^d =
			\left(\int_{-1}^1 \frac{e^{-u^2/2}}{(2\pi)^{1/2}}
			\, du\right)^d f_\e(\s^2)\\
		&:= c_2 f_\e(\s^2).
	\end{split}\end{equation}
	The lemma follows with $c:= \min(c_1\, ,c_2)$.
\end{proof}

Next we find bounds for $F_\e$ in terms of the function
$U_{(d-2)/2}$ that is defined in \eqref{eq:U}.

\begin{lem}\label{lem:f}
	There exists $c>1$ such that
	such that for all $0\le y\le 2$ and $\e>0$,
	\begin{equation}
		F_\e(y) \le
		c \e^d U_{(d-2)/2} (y).
	\end{equation}
	In addition, for all $y\ge \e^2$,
	\begin{equation}
		F_\e(y) \ge \frac{\e^d}{c} U_{(d-2)/2}(y).
	\end{equation}
\end{lem}

\begin{proof}
	Evidently,
	\begin{equation}\label{eq:ff}
		F_\e(y) = \int_0^1 f_\e(x+y)\, dx 
		\le\e^d \int_0^1 \frac{dx}{(x+y)^{d/2}}
		=\e^d\int_y^{1+y} \frac{dx}{x^{d/2}},
	\end{equation}
	and this is an equality when $y \ge \e^2$.
	The remainder of the proof is a direct computation.
\end{proof}

 As regards the functions $G_\e$,
 we first note that 
 \begin{equation}\label{eq:Asreg}
	 G_\e(x) = \mathop{\iint}_{[0,1]^2} f_\e(x+|y|)\, dy.
\end{equation}
The following captures a more useful property of
$G_\e$.

\begin{lem}\label{lem:G}
	There exists $c>1$ such that 
	for all $0<x\le 2$ and $\e>0$,
	\begin{equation}
		G_\e(x) \le c \e^d U_{(d-4)/2}(x).
	\end{equation}
	If, in addition, $x\ge \e^2$ then
	\begin{equation}
		G_\e(x) \ge \frac{\e^d}{c} U_{(d-4)/2}(x).
	\end{equation}
\end{lem}

Lemma \ref{lem:G} follows from Lemma
\ref{lem:f} and one or two
elementary and direct computations.

We conclude this section with a final technical
lemma.

\begin{lem}\label{lem:GScale}
	For all $x,\e>0$,
	\begin{equation}
		G_\e(x) \ge \frac12 \int_0^2 F_\e(x+y)\, dy.
	\end{equation}
\end{lem}

\begin{proof}
	We change variables to find that
	\begin{equation}
		\int_0^2 F_\e(x+y)\, dy =
		\frac12 \int_0^1 F_\e \left( x+\frac y2\right)\, dy
		\ge
		\frac12 \int_0^1 F_\e (x+y)\, dy,
	\end{equation}
	by monotonicity. This proves the lemma.
\end{proof}

\section{Proof of Theorem \ref{thm:BSMain-Intro}}
	\label{sec:BSMain}

In light of \eqref{eq:Kh:Cap-d/2} it suffices to
to prove that 
\begin{equation}\label{capcap}
	\text{Cap}_{d/2}([0\,,1]\times F)>0 \quad
	\text{iff}\quad \text{Cap}_{(d/2)-1}(F)>0.
\end{equation}
The following harmonic-analytic fact does the job, and a little more; it
must be well known, but we could not find it
in a suitable form in the literature.

Recall that a function $f:\R^n\to[0\,,\infty]$ is of
\emph{strict positive type} if: (i) $f$ is locally integrable
away from $0\in\R^n$;
and (ii) the Fourier transform of $f$ is \emph{strictly} positive.
Corresponding to such a function $f$
we can define a function $\Pi_mf$ [equivalently,
the operator $\Pi_m$] as follows:
\begin{equation}
	(\Pi_mf)(x) := \mathop{\int}_{[0,1]^m}
	f(x+y)\, dy\qquad
	\text{for all }x\in\R^{n-m}.
\end{equation}
It is easy to see that
\begin{equation}\label{eq:Af}
	(\Pi_mf)(x) := \mathop{\iint}_{[0,1]^m\times [0,1]^m}
	f(x+y-z)\, dy\, dz\qquad
	\text{for all }x\in\R^{n-m}.
\end{equation}
This is a direct computation when $m=1$; the general
case is proved by induction.
Then, we have

\begin{thm}[Projection theorem for capacities]\label{thm:proj}
	Let $n>1$ be an integer, and suppose that
	$f:\R^n\to[0\,,\infty]$ is of strict positive
	type and continuous on $\R^n\setminus\{0\}$.
	Then, for all all integers $1\le m<n$ and
	compact sets $F\subset\R^{n-m}$,
	\begin{equation}
		\text{\rm Cap}_f \left(
		[0\,,1]^m \times F\right)=
		\text{\rm Cap}_{\Pi_mf} (F).
	\end{equation}
\end{thm}

The proof is divided into two parts. The first part is easier,
and will be dispensed with first.

\begin{proof}[Proof of Theorem \ref{thm:proj} (The Upper Bound).]
	Let $\lambda_m$ denote the Lebesgue measure on $[0\,,1]^m$,
	normalized to have mass one.
	If $\mu\in\mathcal{P}(F)$, then	evidently,
	\begin{equation}
		I_{\Pi_mf} (\mu) =
		I_f (\lambda_m\times \mu) \ge \inf_{\nu\in\mathcal{P}
		([0,1]^m\times F)} I_f(\nu).
	\end{equation}
	The equality follows from \eqref{eq:Af}
	and the theorem of Fubini--Tonelli.
	But it is clear that $\lambda_m\times \mu
	\in \mathcal{P}([0\,,1]^m\times F)$, 
	whence
	$\text{Cap}_{\Pi_mf}(F) \le \text{Cap}_f([0\,,1]^m\times F)$.
	This completes our proof.
\end{proof}

We need some preliminary developments for the
lower bound. For this portion, we identify
the hypercube $[0\,,1)^m$ with the $m$-dimensional
torus $\mathbf{T}^m$ in the usual way. In particular,
note that $\mathbf{T}^m$ is compact in the resulting
quotient topology. Any probability
measure $\mu$ on $[0\,,1)^m\times F$ can be
identified with a probability measure on $\mathbf{T}^m\times F$
in the usual way. We continue to write the latter measure as
$\mu$ as well. Throughout the remainder of this section,
$f:\R^n\to[0\,,\infty]$ is a fixed function
of strict positive type that is also
continuous on $\R^n\setminus\{0\}$.

\begin{lem}\label{lem:equil}
	Suppose $\mathbf{T}^m\times F$ has positive
	$f$-capacity. Then,
	there exists a probability measure
	$\mathbf{e}_{_{\mathbf{T}^m\times F}}$---the ``equilibrium measure''---on 
	$\mathbf{T}^m\times F$ such that
	\begin{equation}
		I_f(\mathbf{e}_{_{\mathbf{T}^m\times F}})=\left[ \text{\rm
		Cap}_f \left( \mathbf{T}^m\times F \right)
		\right]^{-1}<\infty.
	\end{equation}
\end{lem}

\begin{proof}
	For all $\e>0$ we can find $\mu_\e\in\mathcal{P}(
	\mathbf{T}^m\times F)$ such that	
	\begin{equation}\label{eq:almostcap}
		I_f(\mu_\e) \le \frac{1+\e}{\text{\rm
		Cap}_f \left( \mathbf{T}^m\times F \right)}.
	\end{equation}
	All $\mu_\e$'s are probability measures on the
	same compact set $\mathbf{T}^m\times F$. 
	Choose an arbitrary weak limit 
	$\mu_0\in\mathcal{P}(\mathbf{T}^m\times F)$
	of the sequence $\{\mu_\e\}_{\e>0}$,
	as $\e\to 0$. It follows
	from Fatou's lemma that
	\begin{equation}\begin{split}
		\liminf_{\e\to 0} I_f(\mu_\e) 
			&\ge \liminf_{\eta\to 0}
			\liminf_{\e\to 0} \mathop{\iint}_{
			\{|x-y|\ge\eta\}}
			f(x-y)\, \mu_\e(dx)\, \mu_\e(dy)\\
		&\ge \liminf_{\eta\to 0}\mathop{\iint}_{
			\{|x-y|\ge\eta\}}
			f(x-y)\, \mu_0(dx)\, \mu_0(dy)\\
		& = I_f(\mu_0).
	\end{split}\end{equation}
	Thanks to \eqref{eq:almostcap}, $I_f(\mu_0)$
	is at most equal to the reciprocal of the
	$f$-capacity of $\mathbf{T}^m\times F$.
	On the other hand, the said capacity is bounded
	above by $I_f(\s)$ for all $\s\in\mathcal{P}(
	\mathbf{}T^m\times F)$, whence follows
	the lemma.
\end{proof}

The following establishes the uniqueness of the
equilibrium measure.

\begin{lem}\label{lem:Eunique}
	Suppose $\mathbf{T}^m\times F$ has positive
	$f$-capacity $\chi$. If $I_f(\mu)=I_f(\nu)=1/\chi$
	for some $\mu,\nu\in\mathcal{P}(\mathbf{T}^m\times F)$, then
	$\mu=\nu=\mathbf{e}_{_{\mathbf{T}^m\times F}}$.
\end{lem}

\begin{proof}
	We denote by $\F$ the Fourier transform on any
	and every (locally compact) abelian group $G$;
	$\F$ is normalized as follows: For all group characters
	$\xi$, and all $h\in L^1(G)$,
	\begin{equation}
		(\F h)(\xi) = \int_G (x,\xi) h(x)\, dx,
	\end{equation}
	where $(x,\xi)$ is the usual duality relation between
	$x\in G$ and the character $\xi$, and ``$dx$''
	denotes Haar measure (normalized to be one if
	$G$ is compact; counting measure if $G$ is discrete;
	and mixed in the obvious way, when appropriate).
	Because $f$ is of positive type and continuous
	away from the origin,
	\begin{equation}\label{eq:kahane}
		 I_f(\mu) = \frac{1}{(2\pi)^n} \int_{\mathbf{T}^m
		 \times \R^{n-m}}
		 (\F f)(\xi) \left| (\F\mu)(\xi)\right|^2\, d\xi;
	\end{equation}
	see Kahane \cite[Eq.\ (5), p.\ 134]{Kahane}.

	Using \eqref{eq:kahane} (say) we can extend
	the definition of $I_f(\kappa)$ to all signed
	measures $\kappa$ that have finite absolute mass.
	We note that $I_f(\kappa)$ is real and non-negative,
	but could feasibly be infinite;
	$I_f(\kappa)$ is strictly positive if $\kappa$ is not identically
	equal to the zero measure. The latter follows from the strict
	positivity of $f$.
	
	Let $\rho$ and $\s$ denote two signed
	measures that have finite absolute mass.
	Then, we can define, formally,
	\begin{equation}\label{eq:kahane1}
		I_f(\s,\rho) := 
		\iint \left[ \frac{f(x-y)+f(y-x)}{2}\right]
		\, \s(dx)\, \rho(dy).
	\end{equation}
	This is well-defined if
	$I_f(|\sigma|\,,|\rho|)<\infty$, for instance.
	Evidently, $I_f(\s,\rho)=I_f(\rho,\s)$ and
	$I_f(\s,\s)=I_f(\s)$.
	Finally, by the Cauchy--Schwarz inequality,
	\begin{equation}\label{eq:CBS}
		\left| I_f(\s,\rho) \right| \le 
		I_f(\s)I_f(\rho).
	\end{equation}
	
	Now suppose to the contrary
	that the $\mu$ and $\nu$ of the statement of
	the lemma are distinct. Then,
	by \eqref{eq:kahane},
	\begin{equation}
		0 < I_f\left( \frac{\mu-\nu}{2}\right) =
		\frac{I_f(\mu)+I_f(\nu)-2I_f(\mu,\nu)}{4}
		= \frac{\chi^{-1}- I_f(\mu,\nu)}{2},
	\end{equation}
	where, we recall, $\chi^{-1}=
	I_f(\mathbf{e}_{_{\mathbf{T}^m\times F}})$ 
	denotes the reciprocal
	of the $f$-capacity of $\mathbf{T}^m\times F$.
	Consequently, $I_f(\mu,\nu)$ is 
	strictly less than $I_f(\mathbf{e}_{_{\mathbf{T}^m\times F}})$.
	From this we can deduce that
	\begin{equation}\begin{split}
		I_f\left( \frac{\mu+\nu}{2}\right) &=
			\frac{I_f(\mu)+I_f(\nu)+2I_f(\mu,\nu)}{4}
			= \frac{\chi^{-1}+ I_f(\mu,\nu)}{2}\\
		&< I_f(\mathbf{e}_{_{\mathbf{T}^m\times F}}) \le I_f\left( \frac{\mu+\nu}{2}\right).
	\end{split}\end{equation}
	And this is a contradiction. Therefore, $\mu=\nu$;
	also $\mu$ is equal to $\mathbf{e}_{_{\mathbf{T}^m\times F}}$ because of the
	already-proved uniqueness together with
	Lemma \ref{lem:equil}. 
\end{proof}

\begin{proof}[Proof of Theorem \ref{thm:proj}
	(The Lower Bound).]
	It remains to prove that
	\begin{equation}\label{goal}
		\text{Cap}_{\Pi_mf} (F) \ge \text{Cap}_f
		\left( [0\,,1]^m \times F\right).
	\end{equation}
	We will prove the seemingly-weaker statement
	that
	\begin{equation}\label{goal:cap}
		\text{Cap}_{\Pi_mf} (F) \ge\text{Cap}_f
		\left( \mathbf{T}^m \times F\right).
	\end{equation}
	This is seemingly weaker because
	$ \text{Cap}_f( \mathbf{T}^m \times F)
	=  \text{Cap}_f( [0\,,1)^m \times F)$.
	But, in fact, our proof will reveal that for all $q>1$,
	\begin{equation}
		\text{Cap}_{\Pi_mf} (F) \ge 
		q^{-m} \text{Cap}_f
		\left( [0\,,q)^m \times F\right).
	\end{equation}
	the right-hand side is at least
	$q^{-m}\text{Cap}_f([0\,,1]^m
	\times F)$. Therefore, we can let $q\downarrow
	1$ to derive \eqref{goal}, and therefore
	the theorem.

	Having our ultimate goal
	\eqref{goal:cap} in mind,
	we can assume without loss of generality  that
	\begin{equation}\label{Cap>0}
		\text{Cap}_f\left( \mathbf{T}^m\times F
		\right)>0,
	\end{equation}
	so that $\mathbf{e}_{_{\mathbf{T}^m\times F}}$ exists and is the unique
	minimizer in the definition of
	$\text{Cap}_f(\mathbf{T}^m\times F)$
	(Lemmas \ref{lem:equil} and \ref{lem:Eunique}).
	
	Let us write any $z\in \mathbf{T}^m
	\times \R^{n-m}$ as $z=(z',z'')$, where
	$z'\in\mathbf{T}^m$ and $z''\in\R^{n -m}$.
	
	For all $a,b\in\mathbf{T}^m\times\R^{n-m}$ define
	$\tau_a(b)=a+b$. We emphasize that the first $m$
	coordinates of $\tau_a(b)$ are formed by
	addition in $\mathbf{T}^m$ [i.e., component-wise addition
	mod 1 in $[0\,,1)^m$], whereas the next $n-m$
	coordinates of $\tau_a(b)$ are formed by
	addition in $\R^{n-m}$. In particular, 
	$\tau_a(\mathbf{T}^m\times F) = \mathbf{T}^m\times
	(a'' +F)$.
	
	For all $a\in\mathbf{T}^m\times\R^{n-m}$,
	$\mathbf{e}_{_{\mathbf{T}^m\times F}}\circ \tau_a^{-1}$
	is a probability measure on $\tau_a(\mathbf{T}^m\times
	F)$. Moreover, it is easy to see that
	$\mathbf{e}_{_{\mathbf{T}^m\times F}}$
	and 
	$\mathbf{e}_{_{\mathbf{T}^m\times F}}\circ \tau_a^{-1}$
	have the same $f$-energy. Therefore, whenever
	$a''=0$, $\mathbf{e}_{_{\mathbf{T}^m\times F}}\circ \tau_a^{-1}$
	is a probability measure on $\mathbf{T}^m\times F$
	that minimizes the $f$-capacity of $\mathbf{T}^m\times F$. The
	uniqueness of $\mathbf{e}_{_{\mathbf{T}^m\times F}}$
	proves that 
	\begin{equation}
		\mathbf{e}_{_{\mathbf{T}^m\times F}}
		=\mathbf{e}_{_{\mathbf{T}^m\times F}}\circ \tau_a^{-1}
		\quad\text{whenever }a''=0.
	\end{equation}
	See Lemma \ref{lem:Eunique}.
	Now let $X$ be a random variable with values
	in $\mathbf{T}^m\times F$ such that
	the distribution of $X$ is $\mathbf{e}_{_{\mathbf{T}^m\times F}}$.
	The preceding display implies that for all
	$a'\in\mathbf{T}^m$, the distribution of
	$(X'+a',X'')$ is the same as that of $(X',X'')$. The uniqueness
	of normalized Haar measure $\lambda_m$ then implies that
	$X'$ is distributed as $\lambda_m$. In fact,
	for all Borel sets $A\subset\mathbf{T}^m$ and
	$B\subset \R^{n-m}$,
	\begin{equation}\begin{split}
		\mathbf{e}_{_{\mathbf{T}^m\times F}} (A\times B)
			&= \P\left\{ X'\in A \, ,
			X''\in B \right\}\\
		& = \int_{\mathbf{T}^m}
			\P\left\{ X'\in a'+A\,,
			X''\in B\right\}\, da'\\
		&= \E\left[ \lambda_m(A-X') ~;~ X'' \in B\right]\\
		&= \lambda_m(A)\P\left\{ X''\in B\right\}
			:= \lambda_m(A) \mu(B).
	\end{split}\end{equation}
	Now we compute directly to find that
	\begin{equation}
		\text{Cap}_f \left( \mathbf{T}^m\times F\right)
		= \frac{1}{I_f\left(
		\lambda_m\times\mu\right)} 
		= \frac{1}{I_{\Pi_mf}(\mu)}
		\le \frac{1}{\inf_{\s\in\mathcal{P}
		(F)} I_{\Pi_m f} (\s)}.
	\end{equation}
	This proves \eqref{goal:cap}, and therefore the
	theorem.
\end{proof}

Finally we are ready to present the following:

\begin{proof}[Proof of  Theorem \ref{thm:BSMain-Intro}]
	The function $U_\alpha$ is of strict positive
	type for all $0<\alpha<d$. The easiest way to see this is
	to merely recall the following well-known fact from
	harmonic analysis: In the sense
	of distributions, 
	$\mathcal{F}U_\alpha= c_{d,\alpha} U_{d-\alpha}$
	for a positive and finite constant $c_{d,\alpha}$
	\cite[Lemma 1, p.\ 117]{Stein}. We note also that
	$U_\alpha$ is continuous away from the origin.
	Thus, we can combine \eqref{eq:Kh:Cap-d/2}
	with Theorem \ref{thm:proj} to find that
	\begin{equation}\label{AFP}
		\P\left\{ \mathcal{Z}_d
		\cap F \neq\varnothing\right\}>0\
		\Longleftrightarrow\
		\text{Cap}_{\Pi_1 U_{d/2}}(F)>0.
	\end{equation}
	But for all $x\ge \e^2>0$,
	\begin{equation}
		\left( \Pi_1 U_{d/2} \right) (x) =
		\int_0^1 \frac{dy}{
		|x+y|^{d/2}} = \frac{F_\e(x)}{\e^d}.
	\end{equation}
	Therefore, in accord with Lemmas \ref{lem:f}
	and \ref{lem:GScale}, 
	\begin{equation}
		c_1 U_{(d-2)/2}(x)\le ( \Pi_1 U_{d/2} ) (x)
		\le c_2 U_{(d-2)/2}(x),
	\end{equation}
	for all $\e>0$ and $x\ge 2\e^2$.
	Because $(c_1,c_2)$ does not depend on $\e$, the 
	displayed bounds are valid for all $x>0$, whence it
	follows that
	\begin{equation}\label{eq:AleF}
		\frac{1}{c_2} \text{Cap}_{(d-2)/2}(F) \le
		\text{Cap}_{\Pi_1 U_{d/2} }(F) \le
		\frac{1}{c_1}\text{Cap}_{(d-2)/2}(F).
	\end{equation}
	This and \eqref{AFP} together prove the theorem.
\end{proof}

\section{Proof of Theorem \ref{thm:Double-Intro}}
	\label{sec:Double}
	Let $B^{(1)}$ and $B^{(2)}$ be two independent
	Brownian sheets in $\R^d$, and define
	for all $\mu\in\mathcal{P}(\R_+)$,
\begin{equation}
	J_\e(\mu) := \frac{1}{\e^d}  \mathop{\iint}_{[1,2]^2}
	\int
	\mathbf{1}_{\mathbf{A} (\e\,;s,t)} \, \mu(ds)\, dt,
\end{equation}
where $\mathbf{A} (\e;a,b)$ is the event
\begin{equation}\label{eq:Lambda}
	\mathbf{A} (\e;a,b) :=\left\{
	|B^{(2)}(a\,,b_2)-B^{(1)}(a\,,b_1)|\le 
	\e \right\},
\end{equation}
for all $1\le a,b_1,b_2\le 2$ and $\e>0$.

\begin{lem}\label{lem:EZ}
	We have
	\begin{equation}
		\inf_{0<\e<1}\inf_{\mu\in\mathcal{P}([1,2])}
		\E\left[ J_\e(\mu) \right] \label{eq:EZ}>0.
	\end{equation}
\end{lem}

\begin{proof}
	The distribution of
	$B^{(2)}(s\,,t_2)-B^{(1)}(s\,,t_1)$ has a
	density function that is bounded below,
	uniformly for all $1\le s,t_1,t_2\le 2$.
\end{proof}

Next we present a bound for the second moment of
$J_\e(\mu)$. For technical reasons, we first alter $J_\e(\mu)$
slightly. Henceforth, we define
\begin{equation}
	\hat{J}_\e(\mu) := \frac{1}{\e^d} \mathop{\iint}_{[1,3]^2}
	\int
	\mathbf{1}_{\mathbf{A} (\e;s,t) } \, \mu(ds)\, dt.
\end{equation}

\begin{lem}\label{lem:EZ^2}
	There exists a constant $c$ such that
	for all Borel probability measures $\mu$ on $\R_+$ 
	and all $0<\e<1$,
	\begin{equation}
		\E\left[ \left( \hat{J}_\e(\mu) 
		\right)^2 \right] \le \frac{cI_{G_\e}(\mu)}{\e^d}
		\le c I_{(d-4)/2}(\mu).
	\end{equation}
\end{lem}

\begin{proof} 
	For all $\e>0$, $1\le s,u\le 2$, and
	$t,v\in[1\,,2]\times[3\,,4]$ define
	\begin{equation}
		P_\e (s,u;t,v) := 
		\P\left( \mathbf{A} (\e;s,t)\cap
		\mathbf{A} (\e;u,v)\right).
	\end{equation}
	We claim that there exists a constant $c_1$---independent
	of $(s\,,u\,,t\,,v\,,\e)$---such that
	\begin{equation}\label{eq:Q}
		P_\e(s,u;t,v) \le c_1 \e^d f_\e(|s-u|+|t-v|).
	\end{equation}
	Lemmas 2.3 and 2.4 of Dalang
	and Khoshnevisan \cite{DalangKh} contain
	closely-related, but non-identical, results.
	
	Let us assume \eqref{eq:Q}
	for the time being and prove
	the theorem. We will establish \eqref{eq:Q} subsequently.
	
	Owing to \eqref{eq:Q} and the Fubini--Tonelli theorem,
	\begin{equation}\begin{split}
		\E\left[ \left( \hat{J}_\e(\mu) 
			\right)^2 \right]  &\le
			\frac{c_1}{\e^d} \iint
			\mathop{\iint}_{[1,3]^2
			\times [[1\,,3]]^2} f_\e(|s-u|+|t-v|)\,
			dt\, dv\, \mu(ds)\, \mu(du)\\
		&\le \frac{c}{\e^d}\iint G_\e(s-u)
			\,\mu(ds)\, \mu(du)\\
		&=\frac{c I_{G_\e}(\mu)}{\e^d}.
	\end{split}\end{equation}
	See \eqref{eq:Asreg}. This is the first inequality of
	the lemma. The second follows from the first
	and Lemma \ref{lem:G}. Now we proceed to derive
	\eqref{eq:Q}.
	
	By symmetry, it suffices to estimate 
	$P_\e(s,u;t,v)$ in the case that
	$s\le u$. Now we
	carry out the estimates in two separate cases.\\
	
	\textit{Case 1.} First we consider the case
	$t_1\le v_1$ and $t_2\le v_2$. Define
	$\hat{B}^{(i)}$ to be the white noise that corresponds
	to the sheet $B^{(i)}$ ($i=1,2$). Then, consider
	\begin{equation}\begin{split}
		&H^{(1)}_1 := \hat{B}^{(1)} \left( [0\,,s]\times 
			[0\,,t_1] \right),\quad H^{(1)}_2 := \hat{B}^{(1)} \left( [0\,,s]\times
			[t_1\,,v_1] \right),\\
		&\qquad H^{(1)}_3 := \hat{B}^{(1)} \left( [s\,,u]\times
			[0\,,v_1] \right),\\
		&H^{(2)}_1 := \hat{B}^{(2)} \left( [0\,,s]\times
			[0\,,t_2] \right),\quad H^{(2)}_2 := \hat{B}^{(2)} \left( [0\,,s]\times
			[t_2\,,v_2] \right),\\
		&\qquad H^{(2)}_3 := \hat{B}^{(2)} \left(  [s\,,u]\times
			[0\,,v_2] \right).
	\end{split}\end{equation}
	Then, the $H$'s are all totally independent Gaussian
	random vectors. Moreover, we can find independent
	$d$-vectors $\{\g^{(i)}_j\}_{1\le i\le 2,
	1\le j\le 3}$ of i.i.d.\ standard-normals such that
	\begin{equation}\begin{split}
		H^{(1)}_1 &= (st_1)^{1/2} \g^{(1)}_1,\ 
			H^{(1)}_2 = (s(v_1-t_1))^{1/2} \g^{(1)}_2,\\
		&\qquad H^{(1)}_3 = (v_1(u-s))^{1/2} \g^{(1)}_3,\\
		H^{(2)}_1 &= (st_2)^{1/2} \g^{(2)}_1,\
			H^{(2)}_2 = (s(v_2-t_2))^{1/2} \g^{(2)}_2,\\
		&\qquad H^{(2)}_3 = (v_2(u-s))^{1/2} \g^{(2)}_3.
	\end{split}\end{equation}
	In addition,
	\begin{equation}\begin{split}
		P_\e(s,u;t,v) &=\P\left\{
			\begin{matrix}
				\left| H^{(2)}_1 - H^{(1)}_1 \right|\le \e\\
				\left| H^{(2)}_1 + H^{(2)}_2 + H^{(2)}_3
					- H^{(1)}_1 - H^{(1)}_2 - H^{(1)}_3 \right|
					\le\e
			\end{matrix}
			\right\}\\
		&\le \P\left\{ \left| H^{(2)}_1 - H^{(1)}_1 \right|\le \e\right\}\\
		&\qquad\times \P\left\{ \left| H^{(2)}_2 + H^{(2)}_3
			- H^{(1)}_2 - H^{(1)}_3 \right|
			\le 2\e \right\}.
	\end{split}\end{equation}
	The first term on the right is equal to
	the following:
	\begin{equation}\label{eq:5.9}
		\P\left\{ (s(t_1+t_2))^{1/2} |\g|\le \e\right\}
		\le c_2\e^d,
	\end{equation}
	where $c_2>0$ does not depend on $(s,t,u,v,\e)$;
	see Lemma \ref{lem:g-f}. Also, the second term is
	equal to the following:
	\begin{equation}\label{eq:5.10}\begin{split}
		&\P\left\{ \left( s(v_2-t_2)
			+v_2(u-s) +s(v_1-t_1)+v_1(u-s)\right)^{1/2} |\g|
			\le 2\e \right\}\\
		&\le \P\left\{ \left( |v-t| +(u-s)\right)^{1/2} |\g|
			\le 2\e \right\}\\
		&\le c_3 f_\e(|u-s|+|t-v|),
	\end{split}\end{equation}
	and $c_3>0$ does not depend on $(s\,,t\,,u\,,v\,,\e)$.
	We obtain \eqref{eq:Q} by combining \eqref{eq:5.9} and
	\eqref{eq:5.10}. This completes the proof of Case 1.\\
	
	\textit{Case 2.}\ Now we consider the case
	that $t_2\ge v_2$ and $t_1\le v_1$. We can replace the $H_i^{(j)}$'s
	of Case 1 with the following:
	\begin{equation}\begin{split}
		&H^{(1)}_1 := \hat{B}^{(1)} \left( [0\,,s]\times 
			[0\,,t_1] \right),\quad H^{(1)}_2 := \hat{B}^{(1)} \left( [0\,,s]\times
			[t_1\,,v_1] \right),\\
		&\qquad H^{(1)}_3 := \hat{B}^{(1)} \left( [s\,,u]\times
			[0\,,v_1] \right),\\
		&H^{(2)}_1 := \hat{B}^{(2)} \left( [0\,,s]\times
			[0\,,v_2] \right),\quad
			H^{(2)}_2 := \hat{B}^{(2)} \left( [0\,,s]\times
			[v_2\,,t_2] \right),\\
		&\qquad H^{(2)}_3 := \hat{B}^{(2)} \left(  [s\,,u]\times
			[0\,,v_2] \right).
	\end{split}\end{equation}
	It follows then that 
	\begin{equation}\begin{split}
		P_\e(s,u;t,v) &= \P\left\{
			\begin{matrix}
				\left| H^{(2)}_1 + H^{(2)}_2 - H^{(1)}_1 \right|\le \e\\
				\left| H^{(2)}_1 + H^{(2)}_3 - H^{(1)}_1
				- H^{(1)}_2 - H^{(1)}_3 \right|\le \e
			\end{matrix}
			\right\}.
	\end{split}\end{equation}
	One can check covariances and see that the density function
	of $H^{(2)}_1-H^{(1)}_1$ is bounded above by a constant
	$c_1>0$ that does not depend on $(s\,,t\,,u\,,v\,,\e)$. Therefore,
	\begin{equation}\begin{split}
		P_\e(s,u;t,v) &\le c_1\int_{\R^d} \P\left\{
			\begin{matrix}
				\left| H^{(2)}_1 + z \right|\le \e\\
				\left|  H^{(2)}_3 - H^{(1)}_2 - H^{(1)}_3 + z\right|\le \e
			\end{matrix}
			\right\}\, dz\\
		&= c_1\int_{\{|w|\le \e\}} \P\left\{
			\left|  H^{(2)}_3 - H^{(2)}_1
			+ H^{(1)}_2 - H^{(1)}_3 + w\right|\le \e \right\}\, dw\\
		&\le c_1(2\e)^d \P\left\{
			\left|  H^{(2)}_3 - H^{(2)}_1
			+ H^{(1)}_2 - H^{(1)}_3 \right|\le 2\e \right\}.
	\end{split}\end{equation}
	The component-wise variance of this particular combination
	of $H^{(i)}_j$'s is equal to
	$(u-s)(v_1+v_2) + s(v_1-t_1 + v_2-t_2)
	\ge (u-s)+|t-v|.$
	Whence follows \eqref{eq:Q} in the present case.
	
	Symmetry considerations, together with Cases 1 and 2,
	prove that \eqref{eq:Q} holds for all possible configurations
	of $(s\,,u\,,t\,,v)$. This completes our proof.
\end{proof}

For all $i\in\{1\,,2\}$ and
$s,t\ge 0$, we define $\mathscr{F}^{(i)}_{s,t}$ to be the $\s$-algebra
generated by $\{B^{(i)}(u\,,v)\}_{0\le u\le s,\, 0\le v\le t}$;
as usual, we can assume that $\mathscr{F}^{(i)}$'s
are complete and right-continuous in the partial order
``$\prec$'' described as follows: For all $s,t,u,v\ge 0$,
$(s\,,t)\prec (u\,,v)$ iff $s\le u$ and $t\le v$. Based on
$\mathscr{F}^{(1)}$ and $\mathscr{F}^{(2)}$, we define
\begin{equation}
	\mathscr{F}_{s;t,v} := \mathscr{F}^{(1)}_{s,t} 
	\vee \mathscr{F}^{(2)}_{s,v}
	\qquad\text{for all } s,t,v\ge 0.
\end{equation}

The following proves that Cairoli's maximal $L^2$-inequality
holds with respect to the family of $\mathscr{F}_{s;t,v}$'s.

\begin{lem}\label{lem:cairoli}
	Choose and fix a number $p>1$.
	Then for all almost surely non-negative random variables
	$Y\in \mathcal{L}^p := L^p(\Omega,\vee_{s,t,v\ge 0}\mathscr{F}_{s;t,v},\P)$,
	\begin{equation}
		\left\| \sup_{s,t,v \in \mathbf{Q}_+}
		\E\left[ Y\,\left|\, \mathscr{F}_{s;t,v}
		\right.\right] \right\|_{\mathcal{L}^p}
		\le \left( \frac{p}{p-1}
		\right)^3 \|Y\|_{\mathcal{L}^p}.
	\end{equation}
\end{lem}

\begin{proof}
	We propose to prove that
	for all $s,s',t,t',v,v'\ge 0$,
	and all bounded random variables
	$Y$ that are $\mathscr{F}_{s';t',v'}$-measurable,
	\begin{equation}\label{eq:comm}
		\E\left[ Y\, \left|\, \mathscr{F}_{s;t,v} \right.\right]
		= \E\left[ Y\, \left|\, \mathscr{F}_{s\wedge s'; t\wedge t',
		v\wedge v'} \right.\right] \quad\text{a.s.}
	\end{equation}
	This proves that the three-parameter filtration
	$\{\mathscr{F}_{s;t,v}\}_{s,t,v\in\mathbf{Q}_+}$ is 
	\emph{commuting}
	in the sense of Khoshnevisan \cite[p.\ 35]{Kh}.
	Corollary 3.5.1 of the same reference
	\cite[p.\ 37]{Kh} would then finish
	our proof. 
	
	By a density argument, it suffices to demonstrate
	\eqref{eq:comm} in the case that
	$Y= Y_1Y_2$, where $Y_1$ and $Y_2$ are bounded,
	and measurable with respect to
	$\mathscr{F}^{(1)}_{s',t'}$ and $\mathscr{F}^{(2)}_{s',v'}$, respectively.
	But in this case, independence implies that almost surely,
	\begin{equation}\label{eq:comm2}
		\E\left[ Y\, \left|\, \mathscr{F}_{s;t,v} \right.\right]
		= \E\left[ Y_1\, \left|\, \mathscr{F}^{(1)}_{s,t} \right.\right] 
		\E\left[ Y_2\, \left|\, \mathscr{F}^{(2)}_{s,v} \right.\right].
	\end{equation}
	By the Cairoli--Walsh commutation
	theorem \cite[Theorem 2.4.1, p.\ 237]{Kh},
	$\mathscr{F}^{(1)}$ and $\mathscr{F}^{(2)}$ are each two-parameter,
	commuting filtrations.
	Theorem 3.4.1 of Khoshnevisan \cite[p.\ 36]{Kh}
	implies that almost surely,
	\begin{equation}\begin{split}
		\E\left[ Y_1\, \left|\, \mathscr{F}^{(1)}_{s,t} \right.\right] 
			 & =  \E\left[ Y_1\, \left|\, \mathscr{F}^{(1)}_{s\wedge s',t
			\wedge t'} \right.\right],\\
		\E\left[ Y_2\, \left|\, \mathscr{F}^{(2)}_{s,v} \right.\right]
			&= 
			\E\left[ Y_2\, \left|\, \mathscr{F}^{(2)}_{s\wedge s',v
			\wedge v'} \right.\right].
	\end{split}\end{equation}
	Plug this into \eqref{eq:comm2} to obtain
	\eqref{eq:comm} in the case that $Y$ has the
	special form $Y_1Y_2$, as described above. The general
	form of \eqref{eq:comm} follows from the mentioned
	special case and density.
\end{proof}

\begin{lem}\label{lem:continuity}
	Choose and fix a number $p>1$.
	Then for all almost surely non-negative random variables
	$Y\in \mathcal{L}^p := L^p(\Omega,\vee_{s,t,v\ge 0}\mathscr{F}_{s;t,v},\P)$,
	we can find a continuous modification of the three-parameter process
	$\{ \E[Y\,|\, \mathscr{F}_{s;t,v}] \}_{s,t,v\ge 0}$. Consequently,
	\begin{equation}
		\left\| \sup_{s,t,v \ge 0}
		\E\left[ Y\,\left|\, \mathscr{F}_{s;t,v}
		\right.\right] \right\|_{\mathcal{L}^p}
		\le \left( \frac{p}{p-1}
		\right)^3 \|Y\|_{\mathcal{L}^p}.
	\end{equation}
\end{lem}

\begin{proof}
	First suppose $Y=Y_1Y_2$ where
	$Y_i\in \mathcal{L}^p(\Omega,\vee_{s,t\ge 0}\mathscr{F}^{(i)}_{s,t},\P)$.
	In this case, \eqref{eq:comm2} holds by independence.
	Thanks to Wong and Zakai
	\cite{WongZakai}, each of the two conditional
	expectations on the right-hand
	side of \eqref{eq:comm2} has a representation 
	in terms of continuous, two-parameter and one-parameter
	stochastic integrals. This proves the continuity of
	$(s\,,t\,,v)\mapsto \E[Y\,|\, \mathscr{F}_{s;t,v}]$ in the case
	where $Y$ has the mentioned special form. In the general
	case, we can find $Y^1,Y^2,\ldots$ such that:
	(i) Each $Y^i$ has the mentioned special form; and
	(ii) $\|Y^n-Y\|_{\mathcal{L}^p} \le 2^{-n}$.
	We can write, for all integers $n\ge 1$,
	\begin{equation}
		\left| \E[Y^{n+1}\,|\, \mathscr{F}_{s;t,v} ] -
		\E[Y^n\, |\, \mathscr{F}_{s;t,v} ] \right|
		\le \sum_{k=n}^\infty \left|
		\E[Y^{k+1} -Y^k\,|\, \mathscr{F}_{s;t,v}] \right|.
	\end{equation}
	Take supremum over $s,t,v\in\mathbf{Q}_+$
	and apply Lemma \ref{lem:cairoli} to find that
	\begin{equation}\begin{split}
		&\sum_{n=1}^\infty \left\| 
			\sup_{s,t,v\in\mathbf{Q}_+}
			\left| \E[Y^{n+1}\,|\, \mathscr{F}_{s;t,v} ] -
			\E[Y^n\, |\, \mathscr{F}_{s;t,v} ] \right|
			\right\|_{\mathcal{L}^p}\\
		&\le c\sum_{n=1}^\infty\sum_{k=n}^\infty \left\| Y^{k+1}
			- Y^k\right\|_{\mathcal{L}^p}<\infty.
	\end{split}\end{equation}
	Because each $\E[Y^n\,|\,\mathscr{F}_{s;t,v}]$ is continuous
	in $(s\,,t\,,v)$, $\E[Y\,|\, \mathscr{F}_{s;t,v}]$
	has a continuous modification. The ensuing 
	maximal inequality
	follows from continuity and Lemma \ref{lem:cairoli}.
\end{proof}

\begin{lem}\label{lem:E(Z|F)}
	There exists a constant $c$ such that the following
	holds outside a single null set:
	For all $0<\e<1$, $1\le a,b_1,b_2\le 2$,
	and $\mu\in\mathcal{P}(\R_+)$,
	\begin{equation}\label{eq:E(Z|F)}
		\E\left[\left. \hat{J}_\e(\mu)\, \right|\, \mathscr{F}_{a;b_1,b_2}
		\right] \ge \frac{c}{\e^d} 
		\mathop{\int}_{F\cap [a,2]} 
		G_\e(s-a) \, \mu(ds)
		\cdot \mathbf{1}_{\mathbf{A} (\e/2;a,b)}.
	\end{equation}
\end{lem}

\begin{rem}
	As the proof will show, we may have to redefine
	the left-hand side of \eqref{eq:E(Z|F)} on a null-set
	to make things work seemlessly. The details are
	standard, elementary probability theory and will
	go without further mention.
\end{rem}

\begin{proof}
	Throughout this proof we write
	$\mathscr{E} := \mathscr{E}_{a;b;\e}(\mu) :=
	\E[ \hat{J}_\e(\mu)\,  |\, \mathscr{F}_{a;b_1,b_2} ]$. Evidently,
	\begin{equation}\label{eq:5.20}
		\mathscr{E}
		\ge \frac{1}{\e^d} \int_{b_1}^3 
		\int_{b_2}^3
		\int_{F\cap [a,2]} \P\left(
		\left. \mathbf{A} (\e;s,t) \, \right|\, \mathscr{F}_{a;b_1,b_2} \right)
		\, \mu(ds)\, dt_2\, dt_1.
	\end{equation}
	A white-noise decomposition 
	implies the following: For all $s\ge a$, $t_1\ge b_1$, and
	$t_2\ge b_2$,
	\begin{equation}\begin{split}
		B^{(1)}(s\,,t_1) &= B^{(1)}(a\,,b_1) + b_1^{1/2}
			W^1_1(s-a) + a^{1/2} W^1_2 (t_1-b_1)\\
		&\qquad + V^1(s-a\,,t_1-b_1),\\
		B^{(2)}(s\,,t_2) &= B^2(a\,,b_2) + b_2^{1/2}
			W^2_1(s-a) +a^{1/2} W^2_2 (t_2-b_2)\\
		&\qquad+ V^2(s-a\,,t_2-b_2).
	\end{split}\end{equation}
	Here: the $W^i_j$'s are standard, linear Brownian motions;
	the $V^i$'s are Brownian sheets; and the collection
	$\{W^i_j,V^i,B^i(a\,,b_i)\}_{i,j=1}^2$ is totally
	independent. By appealing to this decomposition in conjunction
	with \eqref{eq:5.20} we can infer
	that the following is a lower bound for $\mathscr{E}$,
	almost surely on the event $\mathbf{A} (\e/2;a,b)$:
	\begin{equation}\begin{split}
		&\frac{1}{\e^d} \int_{b_1}^3\int_{b_2}^3 \mathop{\int}_{
			F\cap [a,2]}\, \mu(ds)\, dt_2\, dt_1\\
		&\quad \times\P\left\{\left|\begin{matrix}
			b_2^{1/2} W^2_1(s-a) + a^{1/2}
			W^2_2(t_2-b_2) + V^2(s-a\,,t_2-b_2)\\
			- b_1^{1/2} W^1_1(s-a) - a^{1/2}
			W^1_2(t_1-b_1) - V^1(s-a\,,t_1-b_1)
			\end{matrix}\right|\le \frac{\e}{2} \right\}\\
		&=\frac{1}{\e^d} \int_{b_1}^3\int_{b_2}^3 \mathop{\int}_{
			F\cap [a,2]}
			\P\left\{ \s\, |\g|\le\frac{\e}{2} \right\}
			\, \mu(ds)\, dt_2\, dt_1.
	\end{split}\end{equation}
	Here, $\g$ is a $d$-vector of i.i.d.\ standard-normals,
	and $\s^2$ is equal to the quantity
	$b_2(s-a)+a(t_2-b_2)
	+(s-a)(t_2-b_2) + b_1(s-a) + a(t_1-b_1) +
	(s-a)(t_1-b_1).$
	The range of possible
	values of $a$ and $b$ is respectively
	$[1\,,2]$ and $[1\,,2]^2$. This means that we can find
	a constant $c>0$---independent
	of $(a\,,b\,,s\,,t)$---such that $\s^2\le c\{
	|s-a|+|t-b|\}$. Apply this bound to the previous display;
	then appeal to Lemma \ref{lem:g-f} to find that
	\eqref{eq:E(Z|F)} holds a.s., but the null-set
	could feasibly depend on $(a\,,b\,,\e)$.
	
	To ensure that the null-set can be chosen independently
	from $(a\,,b\,,\e)$, we first note that the integral
	on the right-hand
	side of \eqref{eq:E(Z|F)} is: (i) Continuous in $\e>0$;
	(ii) independent of $b\in[1\,,2]^2$; and
	(iii) lower semi-continuous in $a\in[1\,,2]$.
	Similarly, $(a\,,b\,,\e)\mapsto \mathbf{1}_{
	\mathbf{A} (\e;a,b)}$ is left-continuous in $\e>0$
	and lower semi-continuous in $(a\,,b)\in[1\,,2]^3$.
	Therefore, it suffices to prove that the left-hand
	side of \eqref{eq:E(Z|F)} is a.s.\ continuous
	in $(a\,,b)\in[1\,,2]^3$, and left-continuous in
	$\e>0$. The left-continuity assertion about $\e>0$
	is evident; continuity in $(a\,,b)$ follows if we could prove
	that for all bounded random variables $Y$,
	$(a\,,b)\mapsto \E\left[ Y\,\left|\,
	\mathscr{F}_{a;b_1,b_2}\right.\right]$
	has an a.s.-continuous  modification.
	But this follows from Lemma \ref{lem:continuity}.
\end{proof}

Next we state and prove a quantitative capacity estimate.

\begin{prop}\label{prop:cap}
	Consider the collection of times of double-points:
	\begin{equation}
		D(\omega) := \left\{ 1\le s\le 2 :\, 
		\inf_{t\in[1,2]^2}
		\left| B^{(2)}(s\,,t_2) - B^{(1)}(s\,,t_1)
		\right| (\omega) =0 \right\}.
	\end{equation}
	Then there exists a constant
	$c>1$ such that
	for all compact, non-random sets $F\subseteq[1\,,2]$,
	\begin{equation}
		\frac 1c \text{\rm Cap}_{(d-4)/2}(F) \le
		\P\left\{ D\cap F\neq\varnothing \right\}
		\le c\text{\rm Cap}_{(d-4)/2}(F).
	\end{equation}
\end{prop}

\begin{proof}
	Define the closed random sets,
	\begin{equation}
		D_\e(\omega) := \left\{ 1\le s\le 2 :\, 
		\inf_{t\in[1,2]^2}
		\left| B^{(2)}(s\,,t_2) - B^{(1)}(s\,,t_1)
		\right| (\omega) \le
		\e  \right\}.
	\end{equation}
	Also, choose and fix a probability measure
	$\mu\in\mathcal{P}(F)$. It is manifest
	that  $D_\e$ intersects $F$ almost surely on the 
	event $\{J_\e(\mu)>0\}$.
	Therefore, we can apply
	the Paley--Zygmund inequality to find that
	\begin{equation}
		\P\left\{ D_\e\cap F \neq\varnothing\right\}
		\ge \frac{\left(\E [ J_\e(\mu) ] \right)^2}{
		\E\left[\left( J_\e(\mu)\right)^2\right]}
		\ge \frac{\left(\E [ J_\e(\mu) ] \right)^2}{
		\E\left[\left( \hat{J}_\e(\mu)\right)^2\right]}.
	\end{equation}
	Let $\e\downarrow 0$ and appeal to compactness
	to find that
	\begin{equation}
		\P\left\{ D\cap F\neq\varnothing\right\}
		\ge \frac{\liminf_{\e\to 0}\left(\E [ J_\e(\mu) ] \right)^2}{
		c I_{(d-4)/2}(\mu)}.
	\end{equation}
	[We have used the second bound of Lemma
	\ref{lem:EZ^2}.] According to Lemma \ref{lem:EZ},
	the numerator is bounded below by a strictly positive
	number that does not depend on $\mu$. Therefore,
	the lower bound of our proposition follows from
	optimizing over all $\mu\in\mathcal{P}(F)$.
	
	In order to derive the upper bound we can assume,
	without any loss in generality, that 
	$\P\{ D_\e\cap F\neq\varnothing\}>0$; for otherwise
	there is nothing to prove. 
	
	For all $0<\e<1$ define
	\begin{equation}
		\tau_\e := \inf\left\{
		s\in F :\, \inf_{t\in[1,2]^2}
		\left| B^{(2)}(s\,,t_2) - B^{(1)}(s\,,t_1)
		\right| \le\e  \right\}.
	\end{equation}
	As usual, $\inf\varnothing :=\infty$. It is easy to see
	that $\tau_\e$ is a stopping time with respect to
	the one-parameter filtration $\{\mathscr{H}_s\}_{s\ge 0}$,
	where
	\begin{equation}
		\mathscr{H}_s := \bigvee_{t,v \ge 0} \mathscr{F}_{s;t,v}
		\qquad \text{for all }s\ge 0.
	\end{equation}
	We note also that there exist $[0\,,\infty]$-valued
	random variables
	$\tau'_\e$ and $\tau''_\e$ such that:
	(i) $\tau'_\e\vee \tau''_\e=\infty$ iff $\tau_\e=\infty$;
	and (ii) almost surely on $\{\tau_\e<\infty\}$,
	\begin{equation}
		\left| B^{(2)} (\tau_\e\,,\tau'_\e) - B^{(1)}
		(\tau_\e\,,\tau''_\e) \right| \le\e.
	\end{equation}
	
	Define
	\begin{equation}
		p_\e:=\P\left\{
		\tau_\e<\infty\right\},
		\quad\text{and}\quad 
		\nu_\e (\bullet) := \P\left\{\left. \tau_\e\in\bullet\,
		\right|\, \tau_\e <\infty\right\}.
	\end{equation}
	We can note that
	\begin{equation}\label{pD}
		\inf_{\e>0} p_\e \ge \P\{ D\cap F\neq\varnothing\},
	\end{equation}
	and this is strictly positive by our earlier assumption. 
	Consequently, $\nu_\e$ is well defined as a 
	classical conditional probability, and
	$\nu_\e\in\mathcal{P}(F)$. Now consider
	the process $\{M^\e\}_{0<\e<1}$ defined
	as follows:
	\begin{equation}
		M^\e_{a;b_1,b_2} :=
		\E\left[\left. \hat{J}_\e(\nu_\e)\, \right|\, 
		\mathscr{F}_{a;b_1,b_2}\right].
	\end{equation}
	Thanks to Lemmas \ref{lem:continuity}
	and \ref{lem:E(Z|F)},
	\begin{equation}\label{eq:CBMS}\begin{split}
		&\E\left[ \sup_{a,b_1,b_2\in\R^3_+}
			\left( M^\e_{a;b_1,b_2}
			\right)^2\right]\\
		&\quad\ge \E\left[ 
			\left( M^\e_{\tau_\e;\tau''_\e,\tau'_\e}
			\right)^2\right]\\
		&\quad \ge \frac{c p_\e}{\e^{2d} } \E\left[\left. \left(
			\int_{F\cap [\tau_\e,2]} G_\e( s-\tau_\e )
			\, \nu_\e(ds) \right)^2
			\,\right|\, \tau_\e <\infty \right]\\
		&\quad\ge \frac{cp_\e }{\e^{2d}}\left(\E\left[\left. 
			\int_{F\cap [\tau_\e,2]} G_\e( s-\tau_\e )
			\, \nu_\e(ds) 
			\,\right|\, \tau_\e <\infty \right]\right)^2.
	\end{split}\end{equation}
	The last line is a consequence of
	the Cauchy--Schwarz inequality. We can bound the
	squared term on the right-hand side as follows:
	\begin{equation}\begin{split}
		&\E\left[\left. 
			\int_{F\cap [\tau_\e,2]} G_\e( s-\tau_\e )
			\, \nu_\e(ds) 
			\,\right|\, \tau_\e <\infty \right]
			\\
		&\quad=  \mathop{\iint}_{\{
			s\in F\cap [u,2]\}} G_\e( s-u )
			\, \nu_\e(ds) \,\nu_\e(du)\\
		&\quad\ge \frac12 \iint G_\e( s-u )\,
			\nu_\e(ds)\, \nu_\e(du)
			=\frac12 I_{G_\e} (\nu_\e).
	\end{split}\end{equation}
	Plug this in \eqref{eq:CBMS},
	and appeal to Lemmas \ref{lem:EZ^2} and \ref{lem:cairoli},
	to find that
	\begin{equation}\begin{split}
		\frac{cp_\e}{4\e^{2d}} \left( I_{G_\e}(\nu_\e) \right)^2
		& \le \E\left[ \sup_{a,b_1,b_2\in\mathbf{Q}_+}
			\left( M^\e_{a;b_1,b_2}
			\right)^2\right]\\
		&\le 2^6 \E\left[
			\left( \hat{J}_\e(\nu_\e) \right)^2\right]
			\le \frac{c}{\e^d} I_{G_\e}(\nu_\e).
	\end{split}\end{equation}
	Solve this, using \eqref{pD}, to find that
	\begin{equation}\label{Dcap}
		\P\{ D\cap F\neq \varnothing\}\le 
		\frac{c}{I_{G_\e}(\nu_\e)}.
	\end{equation}
	Choose and fix a number $\eta > 0$. In accord with 
	Lemma \ref{lem:G},
	\begin{equation}
		I_{G_\e}(\nu_\e) \ge \mathop{\iint}_{\{
		|s-u|\ge\eta\}} U_{(d-4)/2}
		(s-u)\, \nu_\e(ds)\,
		\nu_\e(du),
	\end{equation}
	for all $0< \e< \eta^{1/2}$.
	Recall that $\{\nu_\e\}_{\e>0}$
	is a net of probability measures on $F$. Because
	$F$ is compact, Prohorov's theorem ensures that
	there exists a subsequential weak limit $\nu_0\in\mathcal{P}(F)$
	of $\{\nu_\e\}_{\e>0}$, as $\e\to 0$. Therefore, we can apply
	Fatou's lemma to find that
	\begin{equation}\begin{split}
		\liminf_{\e\to 0} I_{G_\e}(\nu_\e)
			&\ge \lim_{\eta\to 0} \mathop{\iint}_{\{
			|s-u|\ge\eta\}} U_{(d-4)/2}
			(s-u)\, \nu_0 (ds)\, \nu_0 (du)\\
		&= I_{(d-4)/2}(\nu_0). 
	\end{split}\end{equation}
	Together with \eqref{Dcap}, the preceding implies 
	that $\P\{D\cap F\neq\varnothing\}$ is at most
	some constant divided by $I_{(d-4)/2}(\nu_0)$.
	This, in turn, in bounded by a constant multiple
	of $\text{Cap}_{(d-4)/2}(F)$. The proposition follows.
\end{proof}

\begin{proof}[Proof of Theorem \ref{thm:Double-Intro}]
	Let $I$ and $J$ be disjoint, closed intervals in $(0\,,\infty)$
	with the added property that $x<y$ for all $x\in I$ and
	$y\in J$. Define
	\begin{equation}
		\mathcal{D}_d(I,J) :=
		\left\{ s>0 :\ B(s\,,t_1)=B(s\,,t_2)
		\text{ for some $t_1\in I$ and $t_2\in J$} \right\}.
	\end{equation}
	We intend to prove that
	\begin{equation}\label{eq:goal:double}
		\P\{\mathcal{D}_d(I,J)\cap F\neq\varnothing\}>0\
		\Longleftrightarrow\ \mathrm{Cap}_{(d-4)/2}(F)>0.
	\end{equation}
	Evidently, this
	implies Theorem \ref{thm:Double-Intro}. Without loss
	of much generality, we may assume that
	$I=[\tfrac12\,,\tfrac32]$, 
	$J=[\tfrac72\,,\tfrac92]$,
	and $F\subseteq[1\,,2]$. Now consider the random
	fields,
	\begin{equation}\begin{split}
		B^{(2)}( s\,, t ) &:= B(s\,,\tfrac52 +t) - B(s\,,\tfrac52)\\
		B^{(1)}( s\,, t ) &:= B(s\,,\tfrac52 -t) - B(s\,,\tfrac52),
	\end{split}\end{equation}
	for $0\le s,t\le 5/2$. Then two covariance computations
	reveal that the random fields
	$\{B^{(1)}(s\,,\frac52-t)-B(s\,,\tfrac52)\}_{1\le s,t\le 2}$
	and $\{B^{(2)}(s\,,\tfrac52 +t)-B^{(2)}(s\,,\tfrac52)\}_{1\le
	s,t\le 2}$ are independent Brownian sheets. On the other hand,
	the following are easily seen to be equivalent:
	(i) There exists $(s\,,t_1\,,t_2)\in[1\,,2]^3$ such that
	$B^{(1)}(s\,,t_1)=B^{(2)}(s\,,t_2)$; and (ii) There exists
	$(s\,,t_1\,,t_2)\in [1\,,2]\times I\times J$ such that
	$B(s\,,t_1)=B(s\,,t_2)$. Therefore, \eqref{eq:goal:double}
	follows from Proposition \ref{prop:cap}. This completes
	our proof.
\end{proof}

\section{More on Double-Points}\label{sec:moreDouble}

Consider the random sets
\begin{equation}\begin{split}
	\hat{\mathcal{D}}_d &:= \left\{ (s\,,t_1\,,t_2)\in\R^3_+:\
		B(s\,,t_1)=B(s\,,t_2) \right\},\\
	\bar{\mathcal{D}}_d &:= \left\{ (s\,,t_1)\in\R^2_+:\
		B(s\,,t_1)=B(s\,,t_2)\text{ for some }t_2>0\right\}.
\end{split}\end{equation}
The methods of this paper are  not sufficiently delicate
to characterize the polar sets of $\hat{\mathcal{D}}_d$
and $\mathcal{D}_d$. I hasten to add that
I believe such a  characterization is within reach of the
existing technology \cite{Kh99}.
Nonetheless it is not too difficult to prove the following
by appealing solely to the techniques developed here.

\begin{thm}\label{D:cap}
	For all non-random compact sets $E\subset(0\,,\infty)^2$
	and $G\subset(0\,,\infty)^3$,
	\begin{equation}\begin{split}
		\text{\rm Cap}_{d/2}(G)>0
			\ \Longrightarrow\
			\P\left\{ \hat{\mathcal{D}}_d \cap G \neq\varnothing
			\right\}>0\ \Longrightarrow\
			\mathcal{H}_{d/2}(G)>0,\\
		\text{\rm Cap}_{(d-2)/2}(E)>0
			\ \Longrightarrow\
			\P\left\{ \bar{\mathcal{D}}_d \cap E \neq\varnothing
			\right\}>0\ \Longrightarrow\
			\mathcal{H}_{(d-2)/2}(E)>0.
	\end{split}\end{equation}
	where $\mathcal{H}_\alpha$ denotes the $\alpha$-dimensional
	Hausdorff measure [Appendix \ref{app:dimh}].
\end{thm}

\begin{proof}
	Let $B^{(1)}$ and
	$B^{(2)}$ be two independent, two-parameter
	Brownian sheets on $\R^d$. 
	It suffices to prove that there exists a constant
	$c>1$ such that for all non-random
	compact sets $E\subseteq[1\,,2]^2$ and
	$G\subseteq[1\,,2]^3$,
	\begin{equation}\label{goal:T}\begin{split}
		c^{-1} \text{\rm Cap}_{d/2}(G) &\le 
			\P\left\{ \hat{\mathcal{T}}_d \cap G \neq\varnothing
			\right\}\le c\mathcal{H}_{d/2}(G),\\
		c^{-1} \text{\rm Cap}_{(d-2)/2}(E) &\le 
			\P\left\{ \bar{\mathcal{T}}_d \cap E \neq\varnothing
			\right\}\le c\mathcal{H}_{(d-2)/2}(E),
	\end{split}\end{equation}
	where 
	\begin{equation}\begin{split}
		\hat{\mathcal{T}}_d := \left\{ 
			(s\,,t_1,t_2)\in [1\,,2]^3:\ 
			B^{(2)}(s\,,t_2)= B^{(1)}(s\,,t_1) \right\},\\
		\bar{\mathcal{T}}_d := \left\{ 
			(s\,,t_1)\in [1\,,2]^2:\ 
			B^{(2)}(s\,,t_2)= B^{(1)}(s\,,t_1) 
			\text{ for some }t_2>0\right\}.
	\end{split}\end{equation}
	[This sort of reasoning has been employed in the
	proof of Theorem \ref{thm:BSMain-Intro} already;
	we will not repeat the argument here.] We begin by
	deriving the first bound in \eqref{goal:T}.
	
	Recall \eqref{eq:Lambda}.
	Choose and fix $\mu\in\mathcal{P}(G)$, and define
	for all $\e>0$,
	\begin{equation}
		\mathcal{J}_\e(\mu) := \frac{1}{\e^d} \iiint
		\mathbf{1}_{\mathbf{A} (\e;s,t)}
		\, \mu(dsdt_1dt_2).
	\end{equation}
	The proof of Lemma \ref{lem:EZ} shows that
	\begin{equation}\label{EJ}
		\inf_{0<\e<1} \inf_{
		\mu\in\mathcal{P}([1,2]^3)}
		\E\left[ \mathcal{J}_\e(\mu)\right] >0.
	\end{equation}
	Similarly, we can apply \eqref{eq:Q} to find that
	\begin{equation}\begin{split}
		\E\left[\left( \mathcal{J}_\e(\mu)\right)^2\right]
			&\le \frac{c}{\e^d} \iiiint f_\e(|s-u|+|t-v|)\, \mu(dsdt_1dt_2)
			\, \mu(dudv_1dv_2)\\
		&\le c I_{d/2}(\mu).
	\end{split}\end{equation}
	We have used the obvious inequality,
	$f_\e(x)\le \e^d |x|^{-d/2}$. The lower bound
	in \eqref{goal:T} follows from the previous two moment-bounds,
	and the Paley--Zygmund--inequality; we omit the details.
	
	For the proof of the upper bound it is convenient to
	introduce some notation. Define
	\begin{equation}\begin{split}
		\Delta(s;t) &:= B^{(2)}(s\,,t_2) - B^{(1)}(s\,,t_1)
			\quad\text{for all }s,t_1,t_2\ge 0,\\
		\mathcal{U}(x;\e) &:= \left[ x_1\,,x_1+\e \right]
			\times \left[ x_2\,,x_2+\e \right] \times
			\left[ x_3\,,x_3+\e\right]\quad
			\text{for all }x\in\R^3,\, \e>0.
	\end{split}\end{equation}
	Then,
	\begin{equation}
		\P\left\{ \hat{\mathcal{T}}_d \cap 
		\mathcal{U}(x;\e) \neq\varnothing \right\}
		\le \P\left\{ \left| \Delta(x) \right| \le
		\Theta(x;\e) \right\},
	\end{equation}
	where $\Theta(x;\e) := \sup_{y\in\mathcal{U}(x;\e)}
	|\Delta(y)-\Delta(x)|$. The density function of
	$\Delta(x)$ is bounded above, uniformly for all
	$x\in[1\,,2]^3$. Furthermore, $\Delta(x)$ is independent
	of $\Theta(x;\e)$. Therefore, there exists a constant
	$c$ such that uniformly for all $0<\e<1$ and
	$x\in[1\,,2]^3$,
	\begin{equation}\label{UU}
		\P\left\{ \hat{\mathcal{T}}_d \cap 
		\mathcal{U}(x;\e) \neq\varnothing \right\}
		\le c\E\left[ \left( \Theta(x;\e) \right)^d \right]
		\le c\e^{d/2}.
	\end{equation}
	The final inequality holds
	because: (i) Brownian-sheet scaling dictates that
	$\Theta(x;\e)$ has the same law as
	$\e^{d/2}\Theta(x;1)$; and (ii)
	$\Theta(x;1)$ has moments of all order,
	with bounds that do not depend on 
	$x\in[1\,,2]^3$ \cite[Lemma 1.2]{OreyPruitt}.
	
	To prove the upper bound we can assume that
	$\mathcal{H}_{d/2}(G)<\infty$. In this case
	we can find $x_1,x_2,\ldots\in[1\,,2]^3$ and $r_1,r_2,\ldots\in(0\,,1)$
	such that $G\subseteq \cup_{i=1}^\infty \mathcal{U}(x_i;r_i)$
	and $\sum_{i=1}^\infty r_i^{d/2}\le 2\mathcal{H}_{d/2}(G)$.
	Thus, by \eqref{UU},
	\begin{equation}\begin{split}
		\P\left\{ \hat{\mathcal{T}}_d \cap 
			G \neq\varnothing \right\} &\le \sum_{i\ge 1}
			\P\left\{ \hat{\mathcal{T}}_d \cap 
			\mathcal{U}(x_i;r_i) \neq\varnothing \right\}\\
		&\le c\sum_{i\ge 1} r_i^{d/2}
			\le 2c \mathcal{H}_{d/2}(G).
	\end{split}\end{equation}
	This completes our proof of the first bound in
	\eqref{goal:T}. 
	
	In order to prove the lower bound for $\bar{\mathcal{T}}_d$
	note that $\bar{\mathcal{T}}_d$ intersects $E$ if and only
	if $\hat{\mathcal{T}}_d$ intersects $[0\,,1]\times E$.
	In \eqref{capcap} we proved that if $E$ is a one-dimensional,
	compact set, then
	$\text{Cap}_{d/2}([0\,,1]\times E)=
	\text{Cap}_{(d-2)/2}(E)$. A similar proof shows
	the same fact holds in any dimension, whence follows the
	desired lower bound for
	the probability that $\bar{\mathcal{T}}_d$
	intersects $E$.
	
	To conclude, it suffices to prove that
	\begin{equation}
		\mathcal{H}_{d/2}([0\,,1]\times E)>0\
		\Longrightarrow\ \mathcal{H}_{(d-2)/2}
		(E)>0.
	\end{equation}
	But this follows readily from Frostman's lemma
	[Appendix  \ref{app:dimh}]. Indeed, the positivity
	of $\mathcal{H}_{d/2}([0\,,1]\times E)$ is equivalent
	to the existence of $\mu\in\mathcal{P}([0\,,1]\times E)$
	and a constant $c$ such that the $\mu$-measure
	of all balls [in $\R^3$] of radius $r>0$ is at most $cr^{d/2}$.
	Define $\bar\mu(C):=\mu([0\,,1]\times C)$ for all
	Borel sets $C\subseteq\R^2$. Evidently,
	$\bar\mu\in\mathcal{P}(E)$, and a covering 
	argument, together with the Frostman property of
	$\mu$, imply that $\bar\mu$ of all two-dimensional
	balls of radius $r>0$ is at most $cr^{(d/2)-1}$.
	Another application of the Frostman lemma finishes
	the proof.
\end{proof}

\section{Proof of Theorem \ref{thm:trace}}
	\label{sec:trace}

Define for all $s>0$, every $\omega\in\Omega$, and all 
Borel sets $I\subseteq \R_+$,
\begin{equation}
	T_d^I(s) (\omega) := \left\{  t\in I :\ B(s\,,t)
	(\omega) =0 \right\}.
\end{equation}
Equivalently, $T_d^I (s) = B^{-1}\{0\} \cap \left( \{s\}\times(0\,,\infty) \right)
\cap I.$
It suffices to prove that for all closed
intervals $I\subset(0\,,\infty)$,
\begin{equation}\label{goal:trace}
	\dimh T_d^I(s) = 0\quad\text{ for all }s>0
	\quad\text{a.s.}
\end{equation}
[N.B.: The order of the quantifiers!]. This, in turn, proves that
\begin{equation}
	\dimh T_d^{\R_+} (s) = \sup_I \dimh T_d^I(s) =0
	\quad\text{for all }s>0,
\end{equation}
where the supremum is taken over all closed intervals
$I\subset(0\,,\infty)$ with rational end-points.
Theorem \ref{thm:trace} follows suit. Without loss
of much generality, we prove \eqref{eq:trace}
for $I:=[1\,,2]$; the more general case follows
from this after a change of notation. To simplify the
exposition, we write
\begin{equation}
	T_d(s) := T_d^{[1,2]}(s).
\end{equation}

Consider the following events:
\begin{equation}
	\mathbf{G}_k(n) := \left\{  
	\sup_{\substack{1\le s,t\le 2\\
	s\le u\le s+(1/k)\\
	t\le v\le t+(1/k)}}
	|B(u,v) - B(s\,,t)| \le n\left(\frac{\log k}{k}\right)^{1/2} \right\},
\end{equation}
where $k,n\ge 3$ are integers.
We will use the following folklore lemma. A 
generalization is spelled out explicitly in
Lacey \cite[Eq.\ (3.8)]{Lacey}.

\begin{lem}\label{lem:Modulus}
	For all $\gamma>0$ there exists $n_0=n_0(\gamma)$
	such that for all $n,k\ge n_0$,
	\begin{equation}
		\P\left( \mathbf{G}_k(n) \right) \ge  1- n_0 k^{-\gamma}.
	\end{equation}
\end{lem}

Next we mention a second folklore result.

\begin{lem}\label{lem:BM}
	Let $\{W(t)\}_{t\ge 0}$ denote a standard
	Brownian motion in $\R^d$. Then, there exists
	a constant $c$ such that
	for all integers $m\ge 1$ and
	$1\le r_1 \le r_2 \le \ldots \le r_m\le 2$,
	\begin{equation}
		\P\left\{ \max_{1\le i\le m}
		|W(r_i)|\le \e \right\} \le 
		c \e^{d}\prod_{2\le i\le m} 
		\left( \frac{\e}{\left( r_i - r_{i-1}\right)^{1/2}}
		\wedge 1\right)^d.
	\end{equation}
\end{lem}

\begin{proof}
	If $|W(r_i)|\le \e$ for all $i\le m$ then
	$|W(r_1)|\le \e$, and 
	$|W(r_i)-W(r_{i-1})|\le 2\e$ for all 
	$2\le i\le m$. Therefore,
	\begin{equation}
		\P\left\{ \max_{1\le i\le m}
		|W(r_i)|\le \e \right\} \le
		\P\left\{ |W(r_1)|\le \e \right\}
		\prod_{2\le i\le m} \P\left\{ 
		|W(r_i -r_{i-1})|\le 2\e \right\}.
	\end{equation}
	A direct computation yields the lemma from this.
\end{proof}

Now define
\begin{equation}
	I_{i,j}(k) := \left[ 1 + \frac ik \,, 1+ \frac{(i+1)}{k} \right]
	\times \left[ 1+ \frac jk \,, 1+\frac{(j+1)}{k} \right],
\end{equation}
where $i$ and $j$ can each run through $\{ 0,\ldots,k-1\}$, and
$k\ge 1$ is an integer. We say that $I_{i,j}(k)$ is
\emph{good} if $I_{i,j}(k)\cap B^{-1}\{0\}\neq\varnothing$.
With this in mind, we define
\begin{equation}
	N_{i,k} := \sum_{0\le j\le k-1} \mathbf{1}_{\{
	I_{i,j}(k) \text{ is good}\}}
\end{equation}

\begin{lem}\label{lem:NLD}
	Suppose $d\in\{2\,,3\}$. Then,
	for all $\gamma>0$ there exists $\alpha=\alpha(d\,,\gamma)>1$
	large enough that
	\begin{equation}
		\max_{0\le i\le k-1}
		\P\left\{ N_{i,k} \ge \alpha (\log k)^{(8-d)/2}   \right\}
		=O\left( k^{-\gamma} \right),
	\end{equation}
	as $k$ tends to infinity.
\end{lem}

\begin{proof}
	On $\mathbf{G}_k(n)$ we have the set-wise inclusion,
	\begin{equation}
		\left\{ I_{i,j}(k)
		\text{ is good} \right\}\subseteq
		\left\{ \left| B\left( 1+ \frac ik
		~,~ 1+ \frac jk \right) \right|\le n\left(
		\frac{\log k}{k}\right)^{1/2} \right\}.
	\end{equation}
	Therefore, for all integer $p\ge 1$,
	\begin{equation}\begin{split}
		&\E\left[ N_{i,k}^p ~;~ \mathbf{G}_k(n) \right]\\
		&\le \mathop{\sum\cdots\sum}_{0\le j_1 ,
			\,\cdots , j_p \le k-1} \P\left\{
			\max_{1\le \ell \le p}
			\left| B\left( 1+ \frac ik ~,~ 1+\frac{j_\ell}{k} \right)
			\right|\le n\left(
			\frac{\log k}{k}\right)^{1/2} \right\}\\
		&= \mathop{\sum\cdots\sum}_{0\le j_1,
			\,\cdots, j_p \le k-1} \P\left\{
			\max_{1\le \ell \le p}
			\left| \left(1+ \frac ik\right)^{1/2} W\left( 1+\frac{j_\ell}{k} \right)
			\right|\le n\left(
			\frac{\log k}{k}\right)^{1/2} \right\}\\
		&\le p! \mathop{\sum\cdots\sum}_{0\le j_1\le
			\cdots\le j_p \le k-1}
			\P\left\{
			\max_{1\le \ell \le p}
			\left| W\left( 1+\frac{j_\ell}{k} \right)
			\right|\le n\left(
			\frac{\log k}{k}\right)^{1/2} \right\},
	\end{split}\end{equation}
	where $W$ denotes a standard $d$-dimensional Brownian motion.
	Because the latter quantity does not depend on the value of
	$i$, Lemma \ref{lem:BM} shows that
	\begin{equation}\begin{split}
		&\max_{0\le i\le k-1}
			\E\left[ N_{i,k}^p ~;~ \mathbf{G}_k(n) \right]\\
		&\qquad\le c p! n^{pd} \left( \frac{\log k}{k}
			\right)^{d/2} \mathop{\sum\cdots\sum}_{0\le j_1\le
			\cdots\le j_p \le k-1} \
			\prod_{2\le\ell\le p}
			\left( {\frac{\log k}{j_\ell -j_{\ell-1}}} \right)^{d/2},
	\end{split}\end{equation}
	for all $k$ large, where we are interpreting
	$1/0$ as one.
	
	Now first consider the case $d=3$. We recall our
	(somewhat unusual) convention about $1/0$,
	and note that
	\begin{equation}\label{eq:3.11}
		\mathop{\sum\cdots\sum}_{0\le j_1\le
		\cdots\le j_p \le k-1}\
		\prod_{2\le\ell\le p}
		(j_\ell -j_{\ell-1})^{-3/2}
		\le k\left( \sum_{l\ge 0}
		\frac{1}{l^{3/2}} \right)^{p-1}.
	\end{equation}
	Therefore, when $d=3$ we can find a constant
	$c_1$---independent of $(p\,,k)$---such that
	\begin{equation}\label{eq:4.12}\begin{split}
		\max_{0\le i\le k-1}
		\E\left[ N_{i,k}^p ~;~ \mathbf{G}_k(n) \right]
		\le p!\, \frac{(c_1\log k)^{3p/2}}{k^{1/2}}
		\le p!\, (c_1\log k)^{3p/2}.
	\end{split}\end{equation}
	By enlarging $c_1$, if need be, we find that this
	inequality is valid for all $k\ge 1$.
	This proves readily that
	\begin{equation}\label{eq:4.13}\begin{split}
		\max_{0\le i\le k-1}
		\E\left[ \exp\left( 
		\frac{ N_{i,k}}{2(c_1 \log k)^{3/2}}
		\right) ~;~ \mathbf{G}_k(n) \right] \le \sum_{p\ge 0} 2^{-p}=2.
	\end{split}\end{equation}
	Therefore, Chebyshev's inequality implies
	that for all $i,k,p\ge 1$ and $a>0$,
	\begin{equation}\label{eq:4.14}\begin{split}
		\max_{0\le i\le k-1}
		\P\left\{  N_{i,k} \ge 2\gamma 
		c_1^{3/2}(\log k)^{5/2} ~;~
		\mathbf{G}_k(n) \right\} \le 2k^{-\gamma}.
	\end{split}\end{equation}
	Note that  $c_1$ may depend on $n$. But we can choose
	$n$ large enough---once and for all---such that
	the probability of the complement of $\mathbf{G}_k(n)$ is
	at most $nk^{-\gamma}$ (Lemma \ref{lem:Modulus}).
	This proves the lemma in the case that $d=3$.
	
	The case $d=2$ is proved similarly, except
	\eqref{eq:3.11} is replaced by
	\begin{equation}
		\mathop{\sum\cdots\sum}_{0\le j_1\le
		\cdots\le j_p \le k-1}\ \prod_{2\le\ell\le p}
		(j_\ell -j_{\ell-1})^{-1}
		\le k\left( \sum_{0\le l\le k}
		\frac 1l \right)^{p-1} \le k(c_2\log k)^{p-1},
	\end{equation}
	where $c_2$ does not depend on $(k\,,p)$,
	and [as before] $1/0:=1$. 
	Equation \eqref{eq:4.12}, when $d=2$, becomes:
	\begin{equation}\begin{split}
		\max_{0\le i\le k-1}
		\E\left[ N_{i,k}^p ~;~ \mathbf{G}_k(n) \right]
		\le p! (c_2\log k)^p.
	\end{split}\end{equation}
	This forms the $d=2$ version of \eqref{eq:4.13}:
	\begin{equation}\begin{split}
		\max_{0\le i\le k-1}
		\E\left[ \exp\left( 
		\frac{ N_{i,k}}{2c_2 \log k}
		\right) ~;~ \mathbf{G}_k(n) \right] \le 2.
	\end{split}\end{equation}
	Thus, \eqref{eq:4.14}, when $d=2$, becomes
	\begin{equation}\begin{split}
		\max_{0\le i\le k-1}
		\P\left\{  N_{i,k} \ge 2\gamma 
		c_2(\log k)^2 ~;~
		\mathbf{G}_k(n) \right\} \le 2k^{-\gamma}.
	\end{split}\end{equation}
	The result follows from this and Lemma \ref{lem:Modulus}
	after we choose and fix a sufficiently large $n$.
\end{proof}

Estimating $N_{i,k}$ is now a simple matter, as the following
shows.

\begin{lem}\label{lem:N}
	If $d\in\{2\,,3\}$, then with probability one,
	\begin{equation}
		\max_{0\le i\le k-1}
		N_{i,k}= O\left( (\log k)^{(8-d)/2}\right)
		\qquad(k\to\infty).
	\end{equation}
\end{lem}

\begin{proof} 
	By Lemma \ref{lem:NLD}, there exists $\alpha>0$
	so large that for all $k\ge 1$ and $0\le i\le k-1$,
	$\P\{ N_{i,k}\ge \alpha (\log k)^{(8-d)/2}\} \le \alpha k^{-3}.$
	Consequently,
	\begin{equation}
		\P\left\{ \max_{0\le i\le k-1}
		N_{i,k}\ge \alpha 
		(\log k)^{(8-d)/2}
		\right\} \le \alpha k^{-2}.
	\end{equation}
	The lemma follows from this and the Borel--Cantelli lemma.
\end{proof}

We are ready to prove Theorem \ref{thm:trace}. As was mentioned
earlier, it suffices to prove \eqref{goal:trace}, and this follows
from our next result.

\begin{prop}\label{pr:NoDim1}
	Fix $d\in\{2\,,3\}$ and define the measure-function
	\begin{equation}
		\Phi(x):= \left[
		\log_+(1/x) \right]^{-(8-d)/2}.
	\end{equation}
	Then, $\sup_{1\le s\le 2}\mathcal{H}_\Phi ( T_d(s) )<\infty$
	a.s.
\end{prop}
The reason is provided by the following elementary lemma
whose proof is omitted.

\begin{lem}
	Suppose $\varphi$ is a measure function such that
	$\liminf_{x\downarrow 0} x^{-\alpha}\varphi(x)=\infty$
	for some $\alpha>0$. Then, for all Borel sets
	$A\subset\R^n$,
	\begin{equation}
		\mathcal{H}_\varphi(A)<\infty\ \Longrightarrow\
		\mathcal{H}_\alpha(A)<\infty\ \Longrightarrow\
		\dimh  A \le\alpha.
	\end{equation}
\end{lem}

Now we prove Proposition~\ref{pr:NoDim1}.

\begin{proof}[Proof of Proposition~\ref{pr:NoDim1}]
	We can construct a generous cover of $T_d(s)$ as follows:
	For all irrational $s\in [i/k\,,(i+1)/k]$, we cover
	$T_d(s)$ intervals of the form
	\begin{equation}
		\left[ 1+ \frac{j}{k} \,, 1+\frac{(j+1)}{k}\right],
	\end{equation}
	where $j$ can be any integer in $\{0,\ldots,k-1\}$ 
	as long as
	$I_{i,j}(k)$ is good. Therefore, for any measure-function $\varphi$,
	\begin{equation}
		\sup_{\substack{1\le s\le 2:\\
		s\text{ is irrational}}}
		\mathcal{H}^{(1/k)}_\varphi \left( T_d(s) \right)
		\le \varphi(1/k) \max_{0\le i\le k-1}
		N_{i,k}.
	\end{equation}
	Now we choose the measure-function
	$\varphi(x):=\Phi(x)$ and let $k\to\infty$ to find that
	$\mathcal{H}_\Phi (T_d(s))$ is finite,
	uniformly over all
	irrational $s\in[1\,,2]$. The case of rational $s$'s
	is simpler to analyse. Indeed, 
	$T_d(s)=\varnothing$ a.s.\ for all rational $s\in[1\,,2]$. This is
	because 
	$d$-dimensional Brownian motion ($d\in\{2\,,3\}$)
	does not hit zero.
\end{proof}

\begin{rem}\label{rem:d=1}
	The form of Lemma \ref{lem:N}
	changes dramatically when $d=1$.
	Indeed, one can adjust the proof of Lemma
	\ref{lem:N} to find that a.s.,
	\begin{equation}
		\max_{0\le i\le k-1}
		N_{i,k}= O\left( k^{1/2}(\log k)^{3/2}\right)
		\qquad(k\to\infty).
	\end{equation}
	This yields fairly readily that the upper Minkowski
	dimension [written as $\dim_{_{\rm M}}$]
	of $T_1(s)$ is at most $1/2$
	simultaneously for all $s>0$. Let $\dim_{_{\rm P}}$
	denote the packing dimension, and recall \eqref{dim:HPM}.
	Then,  the preceding and the theorem of Penrose \cite{Penrose}
	together prove that  almost surely,
	\begin{equation}\label{penrose:HPM}
		\dimh T_1(s) = \dim_{_{\rm P}} T_1(s) =
		\dim_{_{\rm M}} T_1(s)=\frac12\quad
		\text{for all }s>0.
	\end{equation}
\end{rem}

\section{On Rates of Escape}
	\label{sec:rates}

Throughout this section, we choose and fix a
non-decreasing and measurable function 
$\psi:(0\,,\infty)\to(0\,,\infty)$ such that
$\lim_{t\to\infty} \psi(t)=\infty$.
Define, for all Borel-measurable sets $F\subset\R$,
\begin{equation}
	\Upsilon_F (\psi) := \int_1^\infty \left[ \frac{\K_F
	\left(1/\psi(x) \right)}{
	\left( \psi(x)\right)^{(d-2)/2}} \wedge 1\right]
	\frac{dx}{x}.
\end{equation}

\begin{thm}\label{thm:EscMink}
	If $d\ge 3$, then for all
	non-random, compact sets $F\subset(0\,,\infty)$,
	the following holds with probability one:
	\begin{equation}\label{eq:EscMink}
		\liminf_{t\to\infty} \inf_{s\in F} \left(\frac{
		\psi(t)}{t}\right)^{1/2} |B(s\,,t)| 
		 = \begin{cases}
			0&\text{if }\Upsilon_F(\psi)=\infty,\\
			\infty&\text{otherwise}.
		\end{cases}
	\end{equation}
\end{thm}

\begin{rem}
	Although the infimum over all $s\in E$ is
	generally an uncountable one, measurability
	issues do not arise. Our proof actually shows that
	the event in \eqref{eq:EscMink} is a subset
	of a null set. Thus, we are assuming tacitly that
	the underlying probability space is complete.
	This convention applies to the next theorem
	as well.
\end{rem}

\begin{defn}
	Let $F\subset(0\,,\infty)$ be non-random and compact,
	and $\psi:(0\,,\infty)\to(0\,,\infty)$ measurable and
	non-decreasing. Then we say that
	\emph{$(F,\psi)\in\text{\rm FIN}_{\text{loc}}$} if
	there exists a denumerable
	decomposition $F=\cup_{n\ge 1}
	F_n$ of $F$ in terms of closed intervals $F_1,F_2,\ldots$---all
	with rational end-points---such that $\Upsilon_{F_n}(\psi)<\infty$
	for all $n\ge 1$.
\end{defn}

This brings us to the main theorem of this section. 
Its proof is a little delicate because we have
to get three different estimates, each of which is valid only
on a certain scale. This proof is motivated by the earlier work of the
author with David Levin and Pedro M\'endez \cite{KLM}.

\begin{thm}\label{thm:thm:EscPack}
	If $d\ge 3$, then for all
	non-random, compact sets $F\subset(0\,,\infty)$,
	the following holds with probability one:
	\begin{equation}
		\inf_{s\in F}\liminf_{t\to\infty}  \left( \frac{
		\psi(t)}{t} \right)^{1/2} |B(s\,,t)| =
		\begin{cases}
			0&\text{if $(F,\psi)
				\not\in\text{\rm FIN}_{\text{loc}}$},\\
			\infty&\text{otherwise}.
		\end{cases}
	\end{equation}
\end{thm}

The key estimate, implicitly referred to earlier, is the following.

\begin{thm}\label{thm:DK}
	If $d\ge 3$ then
	there exists a constant $c$ such that for all non-random
	compact sets $F\subseteq[1\,,2]$ and $0<\e<1$,
	\begin{equation}
		\frac 1c \left[ \e^{d-2} \K_F(\e^2) \wedge 1\right]
		\le \P\left\{ \inf_{s\in F}
		\inf_{1\le t\le 2} |B(s\,,t)|\le\e \right\} \le
		c \left[ \e^{d-2} \K_F(\e^2) \wedge 1 \right].
	\end{equation}
\end{thm}

Let us mention also the the next result
without proof; it follows upon combining
Theorems 4.1 and 4.2 of our collaborative
effort with Robert Dalang \cite{DalangKh},
together with Brownian scaling:

\begin{lem}\label{lem:DK}
	If $d\ge 3$, then there
	exists $c$ such that for all $1\le a<b\le 2$, $0<\e<1$,
	and $n\ge 1$ such that $(b-a)\ge c\e^2$,
	\begin{equation}
		\frac 1c (b-a)^{(d-2)/2}
		\le \P\left\{ \inf_{\substack{a\le s\le b\\
		1\le t\le 2}} |B(s\,,t)|\le\e \right\} \le
		c (b-a)^{(d-2)/2}.
	\end{equation}
\end{lem}

\begin{rem}
	Dalang and Khoshnevisan \cite{DalangKh} state
	this explicitly for $d\in\{3\,,4\}$. However, the
	key estimates are their Lemmas 2.1 and 2.6,
	and they require only that $d>2$.
\end{rem}

\begin{proof}[Proof of Theorem \ref{thm:DK} (The Upper Bound)]
	Fix $n\ge 1$. Define $I_j := [j/n\,,(j+1)/n)$,
	and let $\chi_j=1$ if $I_j\cap F\neq\varnothing$
	and $\chi_j=0$ otherwise. Then in accord with Lemma \ref{lem:DK},
	\begin{equation}\begin{split}
		&\P\left\{ \inf_{s\in F} \inf_{1\le t\le 2}
			|B(s\,,t)|\le \frac{1}{(cn)^{1/2}}\right\}\\
		&\le \sum_{n\le j\le 2n-1}
			\P\left\{ \inf_{s\in I_j} \inf_{1\le t\le 2}
			|B(s\,,t)|\le \frac{1}{(cn)^{1/2}}\right\} \chi_j\\
		& \le c n^{-(d-2)/2} \M_n(F).
	\end{split}\end{equation}
	This, in turn, is bounded above
	by $c n^{-(d-2)/2} \K_F(1/n)$;
	see \eqref{eq:KMK}.
	The lemma follows in the case that 
	$\e = (cn)^{-1/2}$. The general case follows from a
	\emph{monotonicity argument}, which we rehash
	(once) for the sake of completeness.
	
	Suppose $(c(n+1))^{-1/2} \le \e \le (cn)^{-1/2}$.
	Then,
	\begin{equation}\begin{split}
		\P\left\{ \inf_{s\in F} \inf_{1\le t\le 2}
			|B(s\,,t)|\le \e\right\}
			&\le \P\left\{ \inf_{s\in F} \inf_{1\le t\le 2}
			|B(s\,,t)|\le \frac{1}{(cn)^{1/2}}\right\}\\
		&\le c n^{-(d-2)/2} \K_F(1/n)\\
		&\le c \e^{d-2} \K_F(c\e^2) .
	\end{split}\end{equation}
	Equation \eqref{eq:KK} implies that
	$\K_F(c\e^2)= O(\K_F(\e^2))$ as $\e\to 0$,
	and finishes our proof
	of the upper bound. 
\end{proof}

Before proving the lower bound, let us mention a worthwhile
heuristic argument. If, in Lemma \ref{lem:DK}, the condition
``$(b-a)\ge c\e^2$'' is replaced by $(b-a)\ll \e^2$, then
the bounds both change to $\e^{d-2}$. This is the probability
that a single Brownian motion hits $\mathcal{B}(0;\e)$ some time
during $[1\,,2]$; compare with Lemma \ref{hit:BM}. 
This suggests that the ``correlation length'' among
the slices is of order $\e^2$. That is, slices that are within $\e^2$ of
one another behave much the same; those that are further apart
than $\e^2$ are nearly independent.  We use
our next result in  order to actually 
prove the latter heuristic.

\begin{prop}\label{pr:2BM}
	If $d\ge 3$ then there exists a constant $c$ such that
	for all  $1\le s,u\le 2$ and $0<\e<1$, if $|u-s|\ge\e^2$ then
	\begin{equation}
		\P\left\{ \inf_{1\le t\le 2} |B(s\,,t)|\le \e ~,~
		\inf_{1\le v\le 2} |B(u\,,v)|\le \e \right\}
		\le c \e^{d-2} |u-s|^{(d-2)/2}.
	\end{equation}
\end{prop}
	
\begin{proof}
	Without loss of generality we may choose and
	fix $2\ge u > s\ge 1$. Now the processes
	$\{B(s\,,t)\}_{t\ge 0}$ and
	$\{B(u\,,v)\}_{v\ge 0}$ can be decomposed
	as follows:
	\begin{equation}\begin{split}
		B(s\,,t) = s^{1/2} Z(t),\qquad
		B(u\,,v) = s^{1/2} Z(v) + (u-s)^{1/2} W(v),
	\end{split}\end{equation}
	where $W$ and $Z$
	are independent $d$-dimensional Brownian
	motions. Thus, we are interested in estimating
	the quantity $p_\e$, where
	\begin{equation}\label{eq:p}\begin{split}
		 p_\e & := \P\left\{ \inf_{1\le t\le 2} |Z(t)|
			 \le \frac{\e}{s^{1/2}} ~,~
			 \inf_{1\le v\le 2} \left|  Z(v)+
			 \left( \frac{u-s}{s}\right)^{1/2} W(v) \right|\le \frac{\e}{
			 s^{1/2}} \right\}\\
		&\le \P\left\{ \inf_{1\le t\le 2} |Z(t)|
			 \le \e \ ,
			 \inf_{1\le v\le 2} \left|  Z(v)+
			 (u-s)^{1/2} W(v) \right|\le \e \right\}.
	\end{split}\end{equation}
	The proposition follows from Lemma \ref{hit:2BM}
	in Appendix C below.
\end{proof}

\begin{proof}[Proof of Theorem \ref{thm:DK} (The Lower Bound)]
	We make a discretization argument, once more. 
	Let $n := \K_F(\e^2)$ and find maximal Kolmogorov points
	$s_1<\cdots < s_n$---all in $F$---such that $s_{i+1}-s_i\ge\e^2$
	for all $1\le i<n$. Define
	\begin{equation}
		J_\e(n) := \sum_{1\le i\le n} \mathbf{1}_{\{
		|B(s_i,t)|\le \e \ \text{for some }t\in[1,2]\}}.
	\end{equation}
	According to Lemma \ref{hit:BM},
	\begin{equation}\label{eq:EJ(n)}
		\frac 1c n \e^{d-2}\le \E\left[ J_\e(n) \right] \le c n \e^{d-2}.
	\end{equation}
	On the other hand, the condition $|s_j-s_i|\ge\e^2$ and Proposition
	\ref{pr:2BM} together ensure that
	\begin{equation}
		\E\left[\left( J_\e(n)\right)^2\right] \le \E[J_\e(n)] + c
		\left( \E[J_\e(n)] \right)^2.
	\end{equation}
	Now to prove the lower bound we first assume
	that $n\e^{d-2}\le 1$.
	The previous display implies then that $\E[(J_\e(n))^2]\le c\E[J_\e(n)]$.
	Combine this inequality
	with \eqref{eq:EJ(n)}
	and the Paley--Zygmund inequality to find that
	\begin{equation}
		\P\left\{ \inf_{s\in F}\inf_{1\le t\le 2} |B(s\,,t)|
		\le \e \right\}  \ge \P\left\{ J_\e(n)>0\right\}
		\ge \frac{ \left( \E[J_\e(n)] \right)^2}{\E[(J_\e(n))^2]}
		\ge cn\e^{d-2}.
	\end{equation}
	On the other hand, if $n\e^{d-2}\ge 1$, then the left-hand
	side is bounded away from zero, by a similar bound.
	This is the desired result.
\end{proof}

\begin{lem}\label{lem:DK1}
	Let $d\ge 3$, and $f:[1\,,2]\to \R^d$ be a fixed, non-random,
	measurable function. Then there exists
	a constant $c$ such that for all integers $1\le k\le n$
	\begin{equation}\begin{split}
		&\P\left\{ \inf_{\substack{1\le s\le k/n\\
			1\le t\le 2}} |B(s\,,t) - f(s) |\le\frac{1}{n^{1/2}} 
			\right\}\\
		&\quad \le c \left( k n^{-(d-2)/2} + 
			\sum_{n\le i\le n+k-1} \left(
			\Omega_{i,n}(f) \right)^{d-2} \right),
	\end{split}\end{equation}
	where for all continuous functions $h$,
	\begin{equation}\label{eq:Omega}
		\Omega_{i,n}(h) := \sup_{i/n \le t\le (i+1)/n}
		\left| h(t) - h(i/n) \right|.
	\end{equation}
\end{lem}

\begin{proof}
	Lemma \ref{lem:DK1} holds
	for similar reasons as does Proposition
	\ref{pr:2BM}, but is simpler
	to prove. Indeed, the probability in question is at most
	\begin{equation}
		\sum_{n\le i\le n+k-1}
		\P\left\{ \inf_{i/n \le s\le(i+1)/n}
		|B(s\,,t)-f(s)|\le \frac{1}{n^{1/2}}\right\}.
	\end{equation}
	This, in turn, is less than or equal to
	\begin{equation}
		\sum_{n\le i\le n+k-1}
		\P\left\{ \inf_{1\le t\le 2}
		\left| B(\tfrac{i}{n}\,,t) \right|
		\le \frac{1}{n^{1/2}}
		+ \sup_{1\le t\le 2} \Omega_{i,n}
		(B(\bullet\,,t)) + \Omega_{i,n}(f) \right\}.
	\end{equation}
	By the Markov property, 
	$B((i/n)\,,\bullet)$ is a $d$-dimensional Brownian
	motion that is independent of 
	$\sup_{1\le t\le 2} \Omega_{i,n}(B(\bullet\,,t))$. 
	Standard modulus-of-continuity bounds
	show that the $L^{d-2}(\P)$-norm
	of $\sup_{1\le t\le 2}\Omega_{i,n}(B(\bullet\,,t))$ is at most
	a constant times $n^{-(d-2)/2}$; the details will be explained
	momentarily. Since
	$(i/n)\ge 1$, these observations, in conjunction with Lemma
	\ref{hit:BM} [Appendix C] imply the lemma. It remains to prove
	that there exists a $c$ such that for all $n\ge 1$,
	\begin{equation}\label{supOmega}
		\max_{n\le i\le 2n} \E\left[ \sup_{1\le t\le 2} \left(
		\Omega_{i,n} (B(\bullet\,,t))
		\right)^{d-2} \right] \le c n^{-(d-2)/2}.
	\end{equation}
	Choose and fix $n\ge 1$, $n\le i\le 2n$,
	and $v\in[i/n\,,(i+1)/n]$. Then the process
	$t\mapsto B(v\,,t)-B(i/n\,,t)$ is manifestly a martingale
	with respect to the filtration generated by the
	infinite-dimensional process $t\mapsto B(\bullet\,, t)$.
	Consequently,  $T\mapsto \sup_{1\le t\le T}
	(\Omega_{i,n} (B(\bullet\,,t)))^{d-2}$ is a sub-martingale,
	and \eqref{supOmega} follows from Doob's
	inequality and Brownian-sheet scaling. This completes
	our proof.
\end{proof}

Lemma \ref{lem:DK1}, together with a monotonicity argument,
implies the following.

\begin{lem}\label{lem:DKh}
	Let $d\ge 3$, and $f:[1\,,2]\to \R^d$ be a fixed, non-random,
	measurable function. Then there exists
	a constant $c$ such that for all $1\le a\le 2$
	and $0<\e<1$,
	\begin{equation}\begin{split}
		&\P\left\{ \inf_{a\le s\le a+\e^2}
			\inf_{1\le t\le 3} |B(s\,,t) - f(s) |\le\e
			\right\}\\
		&\quad \le c \left( \e^{d-2} + 
			\sup_{a\le u\le a+\e^2}
			|f(u) - f(a)|^{d-2} \right),
	\end{split}\end{equation}
\end{lem}

\begin{proof}[Proof of Theorem \ref{thm:EscMink}]
	First, assume that $\Upsilon(\psi)<\infty$; this is
	the \emph{first half}.
	
	Define for all $n=0,1,2,\ldots$,
	\begin{equation}\begin{split}
		\psi_n &:= \psi(2^n),\\
		\mathbf{A} _n &:= \left\{ \inf_{s\in F}\inf_{2^n
			\le t\le 2^{n+1}} |B(s\,,t)|
			\le  (2^n/\psi_n )^{1/2} \right\}.
	\end{split}\end{equation}
	We combine Theorem \ref{thm:DK} with the Brownian-sheet
	scaling to deduce the following:
	\begin{equation}\label{eq:P(A_n)}
		\frac 1c \left[
		\psi_n^{-(d-2)/2} \K_F (1/\psi_n) \wedge 1
		\right]
		\le \P(\mathbf{A}_n) \le c 
		\left[ \psi_n^{-(d-2)/2} \K_F(1/\psi_n) \wedge 1\right].
	\end{equation}
	After doing some algebra we find that because
	$\Upsilon_F(\psi)$ is finite, then so is the quantity
	$\sum_{n\ge 1} \P(\mathbf{A}_n)$. By the Borel--Cantelli lemma,
	\begin{equation}
		\liminf_{n\to\infty} \left( \frac{\psi_n}{2^n}\right)^{1/2}
		\inf_{s\in F}\inf_{2^n \le t\le 2^{n+1}}
		|B(s\,,t)| \ge 1\quad\text{a.s.}
	\end{equation}
	If $2^n\le t\le 2^{n+1}$ then $(\psi_n/2^n)^{1/2}\le
	(2\psi(t)/t)^{1/2}$. It follows that almost surely,
	\begin{equation}
		\liminf_{t\to\infty}
		\left( \frac{\psi(t)}{t}\right)^{1/2} 
		\inf_{s\in F} |B(s\,,t)|
		\ge \frac{1}{2^{1/2}}.
	\end{equation}
	But if $\Upsilon_F(\psi)$ is finite then so is
	$\Upsilon_F(r\psi)$, for any $r>0$; see \eqref{eq:KK}.
	Therefore, we can apply the preceding to $r\psi$
	in place of $\psi$, and then let $r\to 0$ to find that
	\begin{equation}
		\Upsilon_F(\psi)<\infty\
		\Longrightarrow\
		\liminf_{t\to\infty}
		\left( \frac{\psi(t)}{t}\right)^{1/2} 
		\inf_{s\in F} |B(s\,,t)|
		=\infty\quad\text{a.s.}
	\end{equation}
	This concludes the proof of the first half.
	
	For the \emph{second half} we assume
	that $\Upsilon_F(\psi)=\infty$. The preceding
	analysis proves that $\sum_{n\ge 1}\P(\mathbf{A}_n)=\infty$.
	According to the Borel--Cantelli lemma, it 
	suffices to prove that
	\begin{equation}\label{goal:BC}
		\limsup_{N\to\infty} \frac{
		\mathop{\sum\sum}_{1\le n<m \le N}
		\P(\mathbf{A}_n \cap \mathbf{A}_m )}{
		\left( \sum_{1\le n\le N} \P(\mathbf{A}_n)
		\right)^2 }<\infty.
	\end{equation}
	
	Define for all integers $n\ge 1$, and all
	$s,t\ge 0$,
	\begin{equation}\begin{split}
		\mathscr{A}_n &:= \text{ the $\s$-algebra
			generated by }\{ B(\bullet\,,v)\}_{0\le v\le 2^n},\\
		\Delta_n (s\,,t) & :=B(s\,,t+2^n) - B(s\,,2^n).
	\end{split}\end{equation}
	The Markov properties of the Brownian sheet
	imply that whenever $m>n\ge 1$:
	(i) $\Delta_m$ is a Brownian sheet that is independent
	of $\mathscr{A}_n$; and (ii)
	$\mathbf{A}_n \in \mathscr{A}_n$. Thus, we apply these
	properties in conjunction with Brownian-sheet scaling
	to find that a.s., $\P ( \mathbf{A}_m\,|\, \mathscr{A}_n )$
	is equal to
	\begin{equation}\begin{split}
		&\P \left( \left. \inf_{s\in F}
			\inf_{2^m-2^n \le t\le 2^{m+1}-2^n}
			\left| \Delta_n(s\,,t) - B(s\,,2^n) \right|
			\le \left(\frac{2^m}{\psi_m}\right)^{1/2} \right| \
			\mathscr{A}_n \right)\\
		&= \P \left( \left. \inf_{1\le t\le (2^{m+1}-2^n)/\alpha}
			\left| \Delta_n(s\,,t) - \frac{B(s\,,2^n)}{
			\alpha^{1/2}} \right|
			\le \left(\frac{2^m}{\alpha \psi_m}\right)^{1/2}
			\ \right| \
			\mathscr{A}_n \right),
	\end{split}\end{equation}
	where $\alpha := 2^m-2^n$. Because $m\ge n+1$,
	$(2^{m+1}-2^n)/\alpha \le 3$
	and $2^m/\alpha\le 2$. Therefore, almost surely,
	\begin{equation}
		\P \left( \mathbf{A}_m\,|\, \mathscr{A}_n \right)
		\le \P \left( \left. \inf_{s\in F}
		\inf_{1\le t\le 3}
		\left| \Delta_n(s\,,t) - \frac{B(s\,,2^n)}{
		\alpha^{1/2}} \right|
		\le \left( \frac{2}{\psi_m} \right)^{1/2}\ \right| \
		\mathscr{A}_n \right).
	\end{equation}
	We can cover $E$ with at most $K:=\M_{[2/\psi_m]}(F)$
	intervals of the form
	$I_i := [i/\ell\,,(i+1)/\ell]$, where $\ell:=[\psi_m/2]$.
	Having done this, a simple bound, together with
	Lemma \ref{lem:DKh} yield the following: With probability one,
	$\P ( \mathbf{A}_m\,|\, \mathscr{A}_n )$ is bounded above by
	\begin{equation}\begin{split}
		&\sum_{1\le i\le K} \P \left( \left. \inf_{s\in I_i}
			\inf_{1\le t\le 3}
			\left| \Delta_n(s\,,t) - \frac{B(s\,,2^n)}{
			\alpha^{1/2}} \right|
			\le \left( \frac{2}{\psi_m}\right)^{1/2} \ \right| \
			\mathscr{A}_n \right)\\
		&\le cK\left( \psi_m^{-(d-2)/2} + 
			\Omega\right),
	\end{split}\end{equation}
	where
	\begin{equation}\begin{split}
		\Omega &:= \alpha^{-(d-2)/2}\max_{1\le i\le K}
			\E\left[ \sup_{s\in I_i}
			\left| B(s\,,2^n) - B(i/\ell\,,2^n) \right|^{d-2} \right]\\
		&= \alpha^{-(d-2)/2} 2^{n(d-2)/2}
			\E\left[ \sup_{0\le s\le 1/\ell}
			| B(s\,,1) |^{d-2} \right]\\
		&= c\alpha^{-(d-2)/2} 2^{n(d-2)/2}\ell^{-(d-2)/2}.
	\end{split}\end{equation}
	Therefore, the bound $2^n/\alpha\le 1$
	implies that $\Omega\le c\ell^{-(d-2)/2}
	\le c \psi_m^{-(d-2)/2}$. On the other hand, 
	by \eqref{eq:KMK} and \eqref{eq:KK}, $K\le \K_F(1/\psi_m)$. Therefore,
	the preceding paragraph and \eqref{eq:P(A_n)}
	together imply that
	$\P(\mathbf{A}_m\,|\,\mathscr{A}_n)\le
	c \P(\mathbf{A}_m)$ a.s., where $c$ does
	not depend on $(n\,,m\,,\omega)$. Therefrom,
	we conclude that
	$\P(\mathbf{A}_m\,|\,\mathbf{A}_n)\le
	c\P(\mathbf{A}_m)$, whence \eqref{goal:BC}.
\end{proof}

We are ready to prove Theorem \ref{thm:thm:EscPack}.

\begin{proof}[Proof of Theorem \ref{thm:thm:EscPack}]
	Suppose, first, that $(F,\psi)\in\text{\rm FIN}_{\text{\it loc}}$.
	According to Theorem \ref{thm:EscMink}, we can 
	write $F=\cup_{n\ge 1} F_n$ a.s., where the $F_n$'s are
	closed intervals with rational end-points, such that 
	\begin{equation}
		\inf_{s\in F_n}
		\liminf_{t\to\infty}  \left( \frac{
		\psi(t)}{t}\right)^{1/2} |B(s\,,t)| = \infty \quad
		\text{for all }n\ge 1.
	\end{equation}
	This proves that a.s.,
	\begin{equation}
		\inf_{s\in F}
		\liminf_{t\to\infty}  \left( \frac{
		\psi(t)}{t}\right)^{1/2} |B(s\,,t)| = \infty,
	\end{equation}
	and this is half of the assertion of the theorem.
	
	Conversely, suppose $(F,\psi)\not\in\text{\rm FIN}_{\text{\it loc}}.$
	Then, given any decomposition $F=\cup_{n\ge 1} F_n$
	in terms of closed, rational intervals $F_1,F_2,\ldots$,
	\begin{equation}\label{nearlythere}
		\liminf_{t\to\infty}  
		\inf_{s\in F_n}\left( 
		\frac{\psi(t)}{t}\right)^{1/2} |B(s\,,t)| = 0 \quad
		\text{for all }n\ge 1.
	\end{equation}
	Define for all $k,n\ge 1$,
	\begin{equation}
		O_{k,n} := \left\{ 
		s>0:\ 
		\inf_{t\ge k} \left[
		\left( \frac{\psi(t)}{t}\right)^{1/2}
		|B(s\,,t)| \right] < \frac 1n \right\}.
	\end{equation}
	Then \eqref{nearlythere} implies that
	every $O_{k,n}$ is relatively open and everywhere
	dense in $F$ a.s. By the Baire category theorem,
	$\cap_{k,n\ge 1} O_{k,n}$ has the same
	properties, and this proves the theorem.
\end{proof}

With Theorem \ref{thm:thm:EscPack} under way,
we can finally derive Theorem \ref{thm:DE} of
the Introduction, and conclude this section.

\begin{proof}[Proof of Theorem \ref{thm:DE}]
	Throughout, define for all $\alpha>0$,
	\begin{equation}
		\psi_\alpha(x) := \left[ \log_+(x) \right]^{2/\alpha}
		\qquad \text{for all }x>0.
	\end{equation}
	
	Note that for any $\psi$, as given
	by Theorem \ref{thm:thm:EscPack}, and for all $\nu>0$,
	\begin{equation}
		\Upsilon_F(\psi)<\infty\quad \text{iff}
		\quad \int_1^\infty \left[ \frac{\K_F
		\left(1/\psi(x) \right)}{
		\left( \psi(x)\right)^{(d-2)/2}} \wedge \nu\right]
		\frac{dx}{x}<\infty.
	\end{equation}
	Therefore,
	\begin{equation}\label{d34}\begin{split}
		&\text{if }\K_F(\e) = O\left( \e^{-(d-2)/2}\right)\
			\ (\e\to 0)\quad\text{then}\\
		&\Upsilon_F(\psi)<\infty\quad\Longleftrightarrow\quad
		\int_1^\infty \frac{\K_F
		\left(1/\psi(x) \right)}{x
		\left( \psi(x)\right)^{(d-2)/2}}\, dx<\infty.
	\end{split}\end{equation}
	
	Suppose $d\ge 4$. Then $K_F(\e)\le c\e^{-1}$, and 
	so by \eqref{d34} and a little calculus,
	\begin{equation}
		\Upsilon_F(\psi_\alpha)<\infty\quad
		\Longleftrightarrow\quad
		\int_1^\infty \frac{\K_F(1/s)}{s^{(d-\alpha)/2}}\,
		ds<\infty.
	\end{equation}
	According to this and \eqref{dim:HPM},
	if $\alpha>d-2-2\dim_{_{\rm M}}F$ is strictly positive,
	then $\Upsilon_F(\psi_\alpha)<\infty$. Theorem
	\ref{thm:EscMink} then implies that, in this case,
	\begin{equation}
		\liminf_{t\to\infty} \inf_{s\in F}
		\frac{(\log t)^{1/\alpha}}{t^{1/2}}
		|B(s\,,t)|=0\quad\text{a.s.}
	\end{equation}
	Similarly, if $0<\alpha<d-2-2\dim_{_{\rm M}}F$, then
	\begin{equation}
		\liminf_{t\to\infty} \inf_{s\in F}
		\frac{(\log t)^{1/\alpha}}{t^{1/2}}
		|B(s\,,t)|=\infty\quad\text{a.s.}
	\end{equation}
	Write $F=\cup_{n\ge 1} F_n$ and ``regularize''
	to find that:
	\begin{enumerate}
		\item If $\alpha > d-2-2\dim_{_{\rm P}}F$ is strictly positive,
			then
			\begin{equation}
				 \inf_{s\in F}\liminf_{t\to\infty}
				\frac{(\log t)^{1/\alpha}}{t^{1/2}}
				|B(s\,,t)|=0\quad\text{a.s.}
			\end{equation}
		\item If $0<\alpha < d-2-2\dim_{_{\rm P}}F$ then
			\begin{equation}
				 \inf_{s\in F}\liminf_{t\to\infty}
				\frac{(\log t)^{1/\alpha}}{t^{1/2}}
				|B(s\,,t)|=\infty\quad\text{a.s.}
			\end{equation}
	\end{enumerate}
	The theorem follows in the case that $d\ge 4$.
	
	When $d=3$, the condition $\dim_{_{\rm M}}F<1/2$
	guarantees that $\K_F(\e)=O(\e^{-1/2})$.
	Now follow through the
	proof of the case $d\ge 4$ to finish. 
\end{proof}

\section{Open Problems}
	\label{sec:open}
	
\subsection{Slices and Zeros}

Theorem \ref{thm:trace} is a metric statement.
Is there a topological counterpart?  The following
is one way to state this formally.

\begin{pbm}\label{pbm1}
	Suppose $d\in\{2\,,3\}$. Is it true that outside a single
	null set, $B^{-1}\{0\}\cap(\{s\}\times(0\,,\infty))$ is
	a finite set for all $s>0$?
\end{pbm}

I conjecture that the answer is ``no.'' In fact, it
is even possible that
there exists a non-trivial measure function $\phi$
such that: (i) $\lim_{r\to 0}\phi(r)=\infty$; and
(ii) $\mathcal{H}_\phi$-measure of 
$B^{-1}\{0\}\cap(\{s\}\times(0\,,\infty))$ is positive
for some $s>0$.

\subsection{Smallness of Double-Points for Slices}
Theorem \ref{D:cap} and a codimension argument
together imply that with probability one,
\begin{equation}\begin{split}
	\dimh \hat{\mathcal{D}}_d &= \left( 3 - \frac d2 \right)_+,
		\quad\text{and}\\
	\dimh \bar{\mathcal{D}}_d &= 2\wedge \left(  3 - \frac d2 \right)_+.
\end{split}\end{equation}
This might suggest that, therefore, none of the slices accrue
any of the dimesion.

\begin{pbm} 
	Define, for all $s\ge 0$,
	\begin{equation}
		\mathcal{Y}_d(s) := \left\{ 
		(t_1,t_2)\in\R^2_+:\ B(s\,,t_1)=B(s\,,t_2)
		\right\}.
	\end{equation}
	Then is it the case that if $d\in\{4\,,5\}$
	then, outside a single null-set,
	$\dimh\mathcal{Y}_d(s)=0$ for all $s\ge 0$?
\end{pbm}
I conjecture that the answer is ``yes.''
Answering this might rely on studying closely the methods of
the literature on ``local non-determinism.'' See,
in particular, Berman \cite{Berman}, Pitt \cite{Pitt},
and the recent deep work of
Xiao \cite{Xiao}. On the other hand,
I believe it should be not too hard to prove
that the answer to the corresponding problem
for $d\le 3$ is ``no,'' due to the existence
of continuous intersection local times \cite{Penrose89}.
[I have not written out a complete
proof in the $d\le 3$ case, mainly because I do not have a
proof, or disproof, in the case that $d\in\{4\,,5\}$. 
This is the more interesting case because
there are no intersection local times.]

Open Problem \ref{pbm1} has the following analogue for
double-points.

\begin{pbm} 
	Let $d\in\{4\,,5\}$. Then is it true that outside
	a single null set, $\mathcal{Y}_d(s)$ is a finite
	set for all $s>0$?
\end{pbm}

The answer to this question is likely to be ``no.'' In fact,
as was conjectured for Open Problem
\ref{pbm1}, here too there might exist
slices that have positive $\mathcal{H}_\phi$-measure
in some gauge $\phi$. If so, then there are in fact
values of $s$ for which $\mathcal{Y}_d(s)$ is
uncountable.

\subsection{Marstrand's Theorem for Projections}
Marstrand \cite{Marstrand} proved that almost every
lower-dimensional orthogonal projection of a Borel set $A$
has the same Hausdorff dimension as $A$.
Theorem \ref{thm:BSMain-Intro} proves that a given
projection (say, onto the $x$-axis) of the zero-set
of Brownian sheet has the same ``Marstrand property.''
I believe that the proof can be adjusted to show that,
in fact, any non-random orthogonal projection of
$B^{-1}\{0\}$ has the same Hausdorff dimension as
$B^{-1}\{0\}$ itself.

\begin{pbm}
	Is there a (random) orthogonal projection such
	that the said projection of $B^{-1}\{0\}$
	has a different Hausdorff dimension than 
	$2-(d/2)$?
\end{pbm}

I believe that the answer is ``no.'' However,
I have no proof nor counter-proof.
Similar questions can be asked about double-points.
I will leave them to the interested reader.

\subsection{Non-Linear SPDEs}

Consider $d$ independent, two-dimensional white noises,
$\dot{B}_1,\ldots,\dot{B}_d$, together with the following
system of $d$ non-interacting stochastic PDEs with
additive noise: For a fixed $T>0$,
\begin{equation}\begin{split}
	\frac{\partial^2 u^i}{\partial t\partial x} (t\,,x) & 
		= \hat{B}_i(t\,,x) + b^i(u(t\,,x)),\\
		u^i(0\,,x) &= u_0(x) \qquad\text{all $-\infty<x<\infty$},\\
	\frac{\partial u^i}{\partial t}(0\,,x) &= u_1(x) \qquad
		\text{all $-\infty<x<\infty$},
\end{split}\end{equation}
where $1\le i\le N$, and $u_0$ and $u_1$ are non-random
and smooth, as well as bounded (say). Then, as long as $b:=(b^1,\ldots,b^d)$
is bounded and Borel-measurable the law of the process
$u:=(u^1,\ldots,u^d)$ is mutually absolutely
continuous with respect to the law of
the two-parameter, $d$-dimensional Brownian sheet $B$.
See Proposition 1.6 of Nualart and Pardoux \cite{NualartPardoux}.
Therefore, the theorems of the preceding sections
apply to the process $u$ equally well.

\begin{pbm}
	Suppose $\s:\R^d\to\R^d$ is a strongly elliptic,
	bounded, $C^\infty$ function.
	Is it the case that the results of the previous
	sections apply to the solution of 
	$(\partial^2 u^i/
	\partial t\partial x)= b^i(u) + \s^i(u)\cdot \hat{B}$
	with reasonable boundary conditions?
\end{pbm}
There is some evidence that the answer is ``yes.''
See Dalang and Nualart \cite{DalangNualart}
where a closely-related problem is solved. 

Finally, we end with an
open-ended question about parabolic SPDEs,
about which we know far less at this point. We will
state things about the additive linear case only.
This case seems to be sufficiently difficult
to analyse at this point in time.

\begin{pbm}
	Consider the following system of
	linear parabolic  SPDE:
	\begin{equation}
		\frac{\partial u^i}{\partial t}(t\,,x) = 
		\frac{\partial^2 u^i}{\partial x^2}(t\,,x) +
		\hat{B}_i(t\,,x),
	\end{equation}
	with reasonable boundary conditions. Is there an analysis of
	the ``slices'' of $u$ along different values of $t$ that is 
	analogous to the results of the present paper?
\end{pbm}

Some results along these lines will appear in
forthcoming work with Robert Dalang and
Eulalia Nualart \cite{DalangKhNualart1,DalangKhNualart2}.

\appendix
\section{Capacity and Dimension}
For the sake of completeness, 
we begin with a brief review of Hausdorff measures. Further
information can be found in Kahane
\cite[Chapter 10]{Kahane},
Khoshnevisan
\cite[Appendices C and D]{Kh}, and Mattila
\cite[Chapter 4]{Mattila}.

\subsection{Capacity}\label{app:cap}

Recall that $\mathcal{P}(F)$ denotes the collection
of all probability measures on the Borel set $F$,
and $|x|$ is the $\ell^1$-norm of the vector $x$.

Let $f:\R^n\to[0\,,\infty]$ be Borel measurable.
Then for all $\mu\in\mathcal{P}(\R^n)$,
the \emph{$f$-energy of $\mu$} is defined by
\begin{equation}
	I_f(\mu) := \iint f(x-y)\, \mu(dx)\, \mu(dy).
\end{equation}
If $F\subset\R^n$ is Borel-measurable, then its
\emph{$f$-capacity} can be defined by
\begin{equation}
	\text{Cap}_f(F) := \left[ \inf_{\mu\in\mathcal{P}(F)}
	I_f(\mu)\right]^{-1},
\end{equation}
where $\inf\varnothing:=\infty$ and $1/\infty:=0$.

Let  $\beta\in\R$ and $x\in\R\setminus\{0\}$ define
\begin{equation}\label{eq:U}
	U_\beta (x) := 
	\begin{cases}
		1,&\text{if $\beta<0$},\\
		\log_+(1/|x|),&\text{if $\beta=0$},\\
		|x|^{-\beta},&\text{if $\beta>0$}.
	\end{cases}
\end{equation}
Also, we define $U_\beta$ at zero by continuously 
extending $U_\beta$
to a $[0\,,\infty]$-valued function on all of $\R$.
Then we write $I_\beta(\mu)$ in place of
$I_{U_\beta}(\mu)$, and $\text{Cap}_\beta(F)$
in place of $\text{Cap}_{U_\beta}(F)$;
$I_\beta(\mu)$ is the \emph{Riesz} [or Bessel--Riesz]
capacity of $\mu$, and $\text{Cap}_\beta$
is [Bessel-] \emph{Riesz capacity} of $F$.

The following is a central property of capacities \cite[p.\ 523]{Kh}.

\begin{taylor}[Taylor \cite{Taylor61}]
	If $F\subset\R^n$ is compact then
	$\mathrm{Cap}_n (F)=0$. Consequently, 
	for all $\beta\ge n$, $\mathrm{Cap}_\beta(F)$ is
	zero also.
\end{taylor}

\subsection{Hausdorff Measures}\label{sec:HM}
A Borel-measurable function $\varphi:\R_+\to[0\,,\infty]$ is said to be a
\emph{measure function} if: (i) $\varphi$ is non-decreasing near
zero; and (ii) $\varphi(2x) = O(\varphi(x))$ as $x\to 0$. Next, 
we choose and
fix a measure function $\varphi$ and a Borel set $A$ in $\R^n$.
For all $r>0$ we define 
\begin{equation}
	\mathcal{H}_\varphi^{(r)}(A) 
	:= \inf \left\{
	\sum_{j\ge 1} \varphi(\delta_j) : \,
	A\subseteq \bigcup_{j\ge 1} \mathcal{B}(x^{(j)};\delta_j),\
	\sup_{j\ge 1}\delta_j\le r,\
	x^{(l)}\in\R^n \right\},
\end{equation}
where $\mathcal{B}(x;r):=\{y\in\R^n: 
|x-y|\le r\}$ is the $\ell^1$-ball of radius
$r>0$ about $x\in\R^n$. 
The \emph{Hausdorff
$\varphi$-measure} $\mathcal{H}_\varphi(A)$ of $A$ can then defined as the
non-increasing limit,
\begin{equation}
	\mathcal{H}_\varphi(A) := \lim_{r\downarrow 0}
	\mathcal{H}_\varphi^{(r)}(A).
\end{equation}
This defines a Borel [outer-] measure on Borel subsets
of $\R^n$.

\subsection{Hausdorff Dimension}\label{app:dimh}
An important special case of $\mathcal{H}_\varphi$
arises when we consider $\varphi(x)=x^\alpha$.
In this case we may write $\mathcal{H}_\alpha$ instead;
this is the \emph{$\alpha$-dimensional Hausdorff
measure}. The \emph{Hausdorff dimension} of $A$
is
\begin{equation}
	\dimh  A := \sup\left\{
	\alpha>0 :\, \mathcal{H}_\alpha(A)>0 \right\}
	=\inf\left\{ \alpha>0 :\, \mathcal{H}_\alpha(A)
	<\infty  \right\}.
\end{equation}
Hausdorff dimension has the following regularity property:
If $A_1,A_2,\ldots$ are Borel sets, then
\begin{equation}
	\dimh \left(
	\bigcup_{i\ge 1} A_i \right)
	= \sup_{i\ge 1}
	\dimh A_i.
\end{equation}
In general, this fails if the union is replaced by an uncountable one.
For instance, consider the example $\R=\cup_{x\in\R}\{x\}$.
The following is another key fact:

\begin{frostman}[Frostman \cite{Frostman}]
	Let $A$ be a compact subset of $\R^n$. Then
	$\mathcal{H}_\alpha(A)>0$ if and only if we
	can find a constant $c$ and a
	$\mu\in\mathcal{P}(A)$ such that
	$\mu(\mathcal{B}(x;r)) \le c r^\alpha$
	for all $r>0$ and $x\in\R^n$.
\end{frostman}
See also Theorem 1 of Kahane \cite[p.\ 130]{Kahane},
Theorem 2.1.1 of
Khoshnevisan \cite[p.\ 517]{Kh}, 
and Theorem 8.8 of Mattila \cite[p.\ 112]{Mattila}.

\section{Entropy and Packing}

The material of this appendix can be found, in expanded
form and with a detailed bibliography, in 
Khoshnevisan et al \cite{KLM}. Throughout, $F\subset\R$
is a Borel-measurable set.

\subsection{Minkowksi Content and Kolmogorov Capacitance}

There are various ways to describe the size of the
set $F$. We have seen already the role
of capacity, Hausdorff measures, and Hausdorff dimension.
Alternatively, we can consider
the rate of growth of the
\emph{Minkowski content} of $F$; this is the function
$\mathbf{N}\ni n\mapsto \M_n(F)$ defined as follows:
\begin{equation}
	\M_n(F) := \#\left\{ i\in\mathbf{Z}:\
	F\cap \left[ \frac in\,,\frac{i+1}{n}\right)\neq
	\varnothing \right\}.
\end{equation}

Also, we can consider the \emph{Kolmogorov
entropy} (known also as ``capacitance''
or ``packing number'')
of $F$; this is the function $(0\,,\infty)\ni \e\mapsto \K_F(\e)$, 
where $\K_E(\e)$ is equal to the maximum number $K$
for which there exists $x_1,\ldots,x_K\in F$
such that $\min_{i\neq j}|x_i-x_j|\ge\e$. Any such sequence
$\{x_i\}_{1\le i\le \K_F(\e)}$ is referred to as a \emph{Kolmogorov
sequence}. 

While $\M_n(F)$ is easier to work with,
$\K_F(\e)$ has the nice property that $\K_F(\e)\ge \K_F(\delta)\ge 1$
whenever $0<\e< \delta$. There are two other properties
that deserve mention. The first is that \cite[Proposition 2.7]{KLM}
\begin{equation}\label{eq:KMK}
	\K_F(1/n) \le \M_n(F) \le 3 \K_F(1/n)\qquad
	\text{for all }n\ge 1.
\end{equation}
The second property
is the following \cite[eq.\ (2.8)]{KLM}:
\begin{equation}\label{eq:KK}
	\K_E(\e) \le 6 \K_F(2\e)\qquad
	\text{for all }\e>0.
\end{equation}

\subsection{Minkowski and Packing Dimension}\label{subsec:dimp}

The (upper) \emph{Minkowski dimension} of $F$ is
the number
\begin{equation}
	\dim_{_{\rm M}} F := \limsup_{n\to\infty} 
	\frac{\log \M_n (F)}{\log n}.
\end{equation}
This is known also as the (upper) ``box dimension'' of $F$,
and gauges the size of $F$. There is a related
\emph{lower Minkowski dimension}; it is defined by
\begin{equation}
	\underline{\dim}_{_{\rm M}} F := \liminf_{n\to\infty}
	\frac{\log \M_n(F)}{\log n}.
\end{equation}

A handicap of the gauge $\dim_{_{\rm M}}$ is that
it assigns the value one to the rationals in $[0\,,1]$;
whereas we often wish to think of
$\mathbf{Q}\cap [0\,,1]$ as a ``zero-dimensional''
set. In such cases, a different notion of dimension
can be used. 

The (upper) \emph{packing dimension} of $F$ is
the ``regularization'' of $\dim_{_{\rm M}}F$ in the
following sense:
\begin{equation}
	\dim_{_{\rm P}} F := \sup\left\{
	\dim_{_{\rm M}} F_k;\ F=\bigcup_{i\ge 1}
	F_i, \text{ $F_i$'s are closed and bounded}\right\}.
\end{equation}
Then it is not hard to see that $\dim_{_{\rm P}}
(\mathbf{Q}\cap [0\,,1])=0$, as desired. Furthermore,
we have the relation,
\begin{equation}\label{dim:HPM}
	\underline{\dim}_{_{\rm M}} F
	\le \dimh F \le \dim_{_{\rm P}} F \le \dim_{_{\rm M}} F.
\end{equation}
These are often equalities; e.g., when $F$ is a
self-similar fractal.
However, there are counter-examples for which
either one, or both, of these inequalities can be strict.
Furthermore, one has \cite[Proposition 2.9]{KLM}
the following integral representations:
\begin{equation}\label{eq:dimint}\begin{split}
	\dim_{_{\rm M}} F & = \inf\left\{ q\in\R: \int_1^\infty
		\K_F(1/s) \, \frac{ds}{s^{1+q}} <\infty
		\right\},\\
	\dim_{_{\rm P}} F & = \inf\left\{ q\in\R:\
	\begin{matrix}
		{}^\exists F_1,F_2,\ldots \text{closed and bounded}\\
		\text{such that } F=\bigcup_{i\ge 1} F_i,\text{ and }\\
		\int_1^\infty s^{-1-q} \K_{F_n}(1/s)\, ds<\infty
		\text{ for all }n\ge 1
	\end{matrix}\right\}.
\end{split}\end{equation}

\section{Some Hitting Estimates for Brownian Motion}

Throughout this section, $X$ and $Y$ denote two
independent, standard Brownian motions in $\R^d$,
where $d\ge 3$.
We will need the following technical lemmas about Brownian motion.
The first lemma is contained in Propositions 1.4.1 and 1.4.3 of 
Khoshnevisan \cite[pp.\ 353 and 355]{Kh}. 

\begin{lem}\label{lem:hit:BM}
	For all $r\in(0\,,1)$,
	\begin{equation}\label{hit:BM}
		\sup_{a\in\R^d}
		\P\left\{\inf_{1\le t\le 3/2} |a+X(t)|\le r\right\}
		\le cr^{d-2} \le c 
		\P\left\{\inf_{1\le t\le 2} |X(t)|\le r\right\}.
	\end{equation}
\end{lem}

We will also need the following variant.

\begin{lem}\label{hit:2BM}
	There exists a
	constant $c$ such that for all 
	$0<r <\rho<1$,
	\begin{equation}\label{eq:hit:2BM}
		\P\left( \left. \inf_{1\le t\le 2}
		\left| \rho Y(t)+X(t) \right| \le r
		\ \right|\ \inf_{1\le s\le 2}|X(s)|\le r \right)
		\le c\rho^{d-2}.
	\end{equation}
\end{lem}

\begin{rem}
	The condition ``$0<\e<\rho<1$'' can be replaced
	with ``$0<\e\le \alpha\rho$'' for any fixed
	finite $\alpha>0$. However, this lemma fails to holds
	for values of $\rho=o(\e)$ as can be seen by simply
	letting $\rho$ tend to zero in the left-hand side of 
	\eqref{eq:hit:2BM}: The left-hand side converges to one
	while the right-hand side converges to zero.
\end{rem}

\begin{proof}
	Define $T:=\inf\{ 1\le t\le 2:\, |X(s)|\le r\}$, where $\inf
	\varnothing :=\infty$, as usual. Then,
	\begin{equation}\begin{split}
		\P_1 & := \P\left( \left. \inf_{T\le t\le 2}  |\rho Y(t)+
			X(t)|\le r\ \right|\ T<\infty\right)\\
		&= \P\left( \left. \inf_{0\le s\le 2-T}
			\left| \rho Y(T+s)+ X(T+s) \right|\le r\ \right|\ 
			T<\infty\right)\\
		&\le \P\left( \left. \inf_{0\le s\le 2-T}
			\left| \rho Y(T+s)+ \hat{X}(s) \right|\le 2r\ \right|\ 
			T<\infty\right),
	\end{split}\end{equation}
	where $\hat{X}(s) := X(T+s)-X(T)$ for all $s\ge 0$.
	By the strong Markov property of $X$,
	\begin{equation}\label{eq:P_1}
		\P_1 \le \sup_{1\le t\le 2} \P\left\{ \inf_{0\le s\le 1}
		\left| \rho Y(t+s)+ X(s) \right|\le 2r\right\}.
	\end{equation}
	In order to estimate this quantity, let us fix an arbitrary
	$t\in[1\,,2]$, and define
	\begin{equation}\begin{split}
		S & := \inf\{ 0\le s\le 1:\ |\rho Y(t+s)+X(s)| \le 2r \},\\ 
		Z &:= \int_0^2 \mathbf{1}_{\{ |\rho Y(t+s)+X(s)|\le 3r\}}\, ds.
	\end{split}\end{equation}
	Then,
	\begin{equation}\begin{split}
		\E[ Z\,|\, S<\infty ] & \ge \E\left[ \left. \int_S^2
			\mathbf{1}_{\{ |\rho Y(t+s)+X(s)|\le 3r\}}\, ds
			\ \right|\ S<\infty \right]\\
		& \ge \E\left[ \left. \int_0^{2-S}
			\mathbf{1}_{\{ |\rho \mathcal{Y}(t+s)+
			\mathcal{X} (s)|\le r\}}\, ds
			\ \right|\ S<\infty \right],
	\end{split}\end{equation}
	where $\mathcal{Y}(u):= Y(u+S)-Y(S)$ and
	$\mathcal{X}(u):= X(u+S)-X(S)$ for all $u\ge 0$.
	The process $u\mapsto \rho Y(t+u)+X(u)$ is
	a L\'evy process, and $S$ is a stopping time
	with respect to the latter process. Therefore, by
	the strong Markov property,
	\begin{equation}\begin{split}
		\E[ Z\,|\, S<\infty ] & \ge \int_0^1
			\P\left\{ |\rho \mathcal{Y}(t+s)+
			\mathcal{X} (s)|\le r\right\} \, ds\\
		&= \int_0^1\P\left\{ \left( \rho^2(t+s)+s\right)^{1/2}
			|\g| \le \e\right\}\, ds\\
		&\ge \int_0^1\P\left\{ \left( \rho^2 t+s\right)^{1/2}
			|\g| \le \e\right\}\, ds,
	\end{split}\end{equation}
	where $\g$ is a $d$-vector of i.i.d.\ standard-normal
	variables. Recall \eqref{eq:f}.
	Thanks to Lemmas \ref{lem:g-f} and \ref{lem:f},
	\begin{equation}
		\inf_{1\le t\le 2} \E[ Z\,|\, S<\infty] \ge
		c \int_0^1 f_\e (\rho^2 +s)\, ds = c F_\e(\rho^2)
		\ge c \e^d \rho^{-(d-2)}.
	\end{equation}
	We have appealed to the condition $\rho>\e$ here.
	Another application of Lemma \ref{lem:g-f}
	yields the following:
	\begin{equation}
		\sup_{1\le t\le 2} \E[Z\,|\, S<\infty] \le
		\frac{\E[Z]}{\P\{S<\infty\}} \le \frac{c\e^d}{\P\{S<\infty\}}.
	\end{equation}
	Recall \eqref{eq:P_1} to find that
	the preceding two displays together imply that
	$\P_1\le c\rho^{d-2}$. Thus, it suffices to prove that
	\begin{equation}\label{eq:P_2}
		\P_2 := \P\left( \left. \inf_{1\le t\le T}  |\rho Y(t)+
		X(t)|\le r\ \right|\ T<\infty\right) \le c\rho^{d-2}.
	\end{equation}
	The estimate on $\P_2$ is derived by using the method
	used to bound $\P_1$; but we apply the latter method
	to the time-inverted Brownian motion $\{ tX(1/t)\}_{t>0}$
	in place of $X$. We omit the numerous, messy details.
\end{proof}

\end{document}